\magnification=\magstep1
\baselineskip=18pt
\font\bbd = msbm10 at 10pt
\font\bbds = msbm8 at 8pt
\font\bbds = msbm7 at 8pt
\font\csc = cmcsc10 at 10pt
\font\small = cmr8 at 8pt
\font\frak = eufm10 at 10pt
\font\fraksmall = eufm8 at 8pt
\font\xmplbx = cmbx10 scaled \magstep1
\outer\def\beginsection#1{\vskip 0pt plus.3\vsize\penalty+50
                          \vskip 0pt plus-.3\vsize\bigskip\vskip\parskip
                          \message{#1}\leftline{\bf #1}
                          \nobreak\smallskip}
\outer\def\proclaim #1 #2\par{\medbreak
      \noindent{\bf#1\enspace}{\sl#2}\par
      \ifdim\lastskip<\medskipamount \removelastskip\penalty55\medskip\fi}
\def\proof {\noindent {\it Proof.} \ \ }

\def\sqr#1#2{{\vcenter{\hrule height.#2pt
              \hbox{\vrule width.#2pt height#1pt \kern#1pt \vrule width.#2pt}
                    \hrule height.#2pt}}}
\def\qe{\mathchoice\sqr54\sqr54\sqr{2.1}3\sqr{1.5}3}

\def\qed{\quad\qe}

\def\smapdown#1{\big\downarrow\rlap{$\vcenter{\hbox{$\scriptstyle#1$}}$}}
\def\C{\mathord{\hbox{\bbd C}}}
\def\c{\mathord{\hbox{\bbds C}}}

\def\F{\mathord{\hbox{\bbd F}}}

\def\V{\mathord{\hbox{\frak V}}}

\def\A{\mathord{\hbox{\frak a}}}

\def\F{\mathord{\hbox{\frak F}}}
\def\E{\mathord{\hbox{\frak E}}}
\def\g{\mathord{\hbox{\fraksmall g}}}
\def\G{\mathord{\hbox{\frak g}}}

\def\H{\mathord{\hbox{\frak h}}}

\def\Z{\mathord{\hbox{\bbd Z}}}
\def\z{\mathord{\hbox{\bbds Z}}}
\def\lsemi{\mathbin{\hbox{\bbd n}}}

\def\setbrackets #1 #2{\bigr\{\,#1\bigm|#2\,\bigr\}}

\def\dd{\mathord{\hbox{$\partial$}}}
\def\dt{\mathord{\hbox{${d\over{d t}}$}}}

\centerline { \xmplbx Extensions of Conformal Modules}
\medskip
\centerline { by }
\medskip
\centerline {Shun-Jen Cheng\footnote{$^1$}{\small partially supported by NSC
grant 86-2115-M-006-012 of the ROC}, Victor G. Kac\footnote{$^{2}$}{\small
partially supported by NSF grant DMS-9622870} and Minoru
Wakimoto$^3$}
\medskip
\vbox{\small\baselineskip=10pt  $^1$Department of Mathematics, National
Cheng-Kung University, Tainan,
Taiwan\hfill\break\indent\phantom{$^1$}chengsj@mail.ncku.edu.tw}
\medskip
\vbox{\small\baselineskip=10pt $^{2}$Department of Mathematics, MIT, Cambridge,
MA 02139, USA\hfill\break\indent\phantom{$^2$}kac@math.mit.edu}
\medskip
\vbox{\small\baselineskip=10pt  $^3$Graduate School of Mathematics, Kyushu
University, Fukuoka 812-81,
Japan\hfill\break\indent\phantom{$^3$}wakimoto@math.kyushu-u.ac.jp}
\medskip

\vskip 0.5in {\noindent\bf Abstract.}  In this paper we classify extensions between irreducible finite conformal modules over the Virasoro algebra, over the current algebras and over their semidirect sum.

\beginsection{1. Introduction}

It was shown in [2] that a finite semisimple conformal algebra is a direct sum of either a current
conformal algebra associated with a finite-dimensional semisimple Lie
algebra, the Virasoro conformal algebra, or their semidirect sum.  The problem of classifying their finite irreducible modules was solved in [1].  However, modules over conformal algebras (and thus conformal modules) are in general not completely reducible, thus in order to understand the representation theory of conformal algebras, one is led to study the extension problem.  This problem is solved in this paper.

In this section we will review a few concepts on conformal algebras, their modules and their relations to Lie algebras of formal distributions that we need in this paper.  For more details the reader is referred to [1], [4], [5].

A {\it formal distribution} (usually called a field by physicists) with
coefficients in a complex vector space $U$ is a generating series of the
form:
$$a(z)=\sum_{n\in\z}a_{[n]}z^{-n-1},$$
where $a_{[n]}\in U$ and $z$ is an indeterminate.

Two formal distributions $a(z)$ and $b(z)$ with coefficients in a Lie
algebra $\G$ are called (mutually) {\it local} if for some $N\in\Z_+$
one has:
$$(z-w)^N[a(z),b(w)]=0.\eqno{(1.1)}$$
Introducing the {\it formal delta function}
$$\delta(z-w)=z^{-1}\sum_{n\in\z}({z\over w})^n,$$
we may write a condition equivalent to (1.1):
$$[a(z),b(w)]=\sum_{j=0}^{N-1}(a_{(j)}b)(w){\dd_w^{j}}{{\delta(z-w)}\over j!},\eqno{(1.2)}$$
for some formal distributions $(a_{(j)}b)(w)$ ([4], Theorem 2.3),
which are uniquely determined by the formula
$$(a_{(j)}b)(w)={\rm Res}_z(z-w)^j[a(z),b(w)].\eqno{(1.3)}$$
Formula (1.3) defines a $\C$-bilinear product $a_{(j)}b$ for each $j\in\Z_+$  
on the space of all formal distributions with coefficients in $\G$.

Note also that the space (over $\C$) of all formal distributions with
coefficients in $\G$ is a (left) module over $\C[\dd_z]$, where the action of  
$\dd_z$ is defined in the obvious way, so that $\dd_z a(z)=\sum_n(\dd
a)_{[n]}z^{-n-1}$, where $(\dd a)_{[n]}=-na_{[n-1]}$.

The Lie algebra $\G$ is called a {\it  Lie algebra of formal
distributions} if there exists a family $\F$ of pairwise local formal
distributions whose coefficients span $\G$. In such a case we say that the
family $\F$ {\it spans} $\G$. We will write $(\G,\F)$ to emphasize the
dependence on $\F$.

The simplest example of a Lie algebra of formal distributions is the
{\it current algebra} $\tilde{\G}$ associated to a finite-dimensional
Lie algebra $\G$:
$$\tilde{\G}=\G\otimes_{\c}\C[t,t^{-1}].$$
It is spanned by the following family of pairwise local formal distributions:
$$a(z)=\sum_{n\in\z}(a\otimes t^n)z^{-n-1},\quad a\in\G.$$
Indeed, it is immediate to check that
$$[a(z),b(w)]=[a,b](w)\delta(z-w).$$

The simplest example beyond current algebras is the (centerless) {\it
Virasoro algebra}, the Lie algebra $\V$ with basis $L_n$ ($n\in\Z$) and
commutation relations
$$[L_m,L_n]=(m-n)L_{m+n}.$$
It is spanned by the local formal distribution $L(z)=\sum_{n\in\z}L_n
z^{-n-2}$, since one has:
$$[L(z),L(w)]=\dd_{w}L(w)\delta(z-w)+2L(w)\dd_{w}\delta(z-w).\eqno{(1.4)}$$
The next important example is the semidirect sum $\V\lsemi \tilde{\G}$ with  
the usual commutation relations between $\V$ and $\tilde{\G}$:
$$[L_m,a\otimes t^n]=-na\otimes t^{m+n},$$
which is equivalent to
$$[L(z),a(w)]=\dd_{w}a(w)\delta(z-w)+a(w)\dd_w\delta(z-w).$$

Given a Lie algebra of formal distributions $(\G,\F)$, we may always
include $\F$ in the minimal family $\bar{\F}$ of pairwise local distributions
which is closed under $\dd$ and all products (1.3) ([4], Section 2.7).  Then  
$\bar{\F}$ is a {\it conformal algebra} with respect to the products (1.3). Its
definition is given below [4]:

A {\it conformal algebra} is a left $\C[\dd]$-module
$R$ with a $\C$-bilinear product $a_{(n)}b$ for each $n\in\Z_+$ such that the  
following axioms hold ($a,b,c\in R; m,n\in\Z_+$):
\item{(C0)} \quad $a_{(n)}b=0$, for $n>>0$,
\item{(C1)} \quad $(\dd a)_{(n)}b=-na_{(n-1)}b,\quad a_{(n)}\dd b=\dd(a_{(n)}b)+n a_{(n-1)}b$,
\item{(C2)} \quad $a_{(n)}b=\sum_{j=0}^{\infty}(-1)^{j+n+1}{{\dd^{j}}\over j!} {(b_{(n+j)}a)}$,
\item{(C3)} \quad $a_{(m)}(b_{(n)}c)=\sum_{j=0}^{\infty}{m\choose
j}(a_{(j)}b)_{(m+n-j)}c+b_{(n)}(a_{(m)}c)$.

Conversely, assuming for simplicity that a conformal algebra $R=\oplus_{i\in  
I}\C[\dd]a^i$ is free as a
$\C[\dd]$-module, we may associate to $R$ a Lie
algebra of formal distributions $(\G(R),\F)$ with basis $a^i_{[m]}$
($i\in I$, $m\in\Z$) and $\F=\{a^i(z)=\sum_na^i_{[n]}z^{-n-1}\}_{i\in I}$
with bracket (cf.~(1.2)):
$$[a^i(z),a^j(w)]=\sum_{k\in\z_+}(a^i_{(k)}a^j)(w){\dd_w^{k}}{{\delta(z-w)}\over k!},\eqno{(1.5)}$$
so that $\bar{\F}=R$. Formula (1.5) is equivalent to the following commutation relations ($m,n\in\Z$; $i,j\in I$):
$$[a^i_{[m]},a^j_{[n]}]=\sum_{k\in\z_+}{m\choose
k}(a^i_{(k)}a^j)_{[m+n-k]}.\eqno{(1.6)}$$

The simplest example of a conformal algebra is the {\it current conformal  
algebra} associated to a finite-dimensional Lie algebra $\G$:
$$R(\tilde{\G})=\C[\dd]\otimes_{\c}\G,$$
with products defined by:
$$a_{(0)}b=[a,b],\quad a_{(j)}b=0,\quad {\rm for\ }j>0,\ a,b\in\G,$$
and the {\it Virasoro conformal algebra} $R(\V)=\C[\dd]\otimes_{\c} L$ with  
products (cf.~(1.4)):
$$L_{(0)}L=\dd L,\quad L_{(1)}L=2L,\quad L_{(j)}L=0,\quad {\rm for\ }j>1.$$
Due to (C1), it suffices to define products on the generators of $R$ over $\C[\dd]$.

Their semidirect sum is $R(\V)\lsemi R(\tilde{\G})=R(\V\lsemi\tilde{\G})$ with additional non-zero
products $L_{(0)}a=\dd a$ and $L_{(1)}a=a$, for $a\in\G$. These examples
are the conformal algebras associated to the Lie algebras of
formal distributions described above.

Let $(\G,\F)$ be a Lie algebra of formal distributions, and let $V$ be  
a $\G$-module.  We say that a formal distribution $a(z)\in\F$ and a formal
distribution $v(z)=\sum_{n\in\z}v_{[n]}z^{-n-1}$ with coefficients in $V$ are  
{\it local} if
$$(z-w)^N a(z)v(w)=0,\quad {\rm for\ some\ }N\in\Z_{+}.\eqno{(1.7)}$$
It follows from [4] {Section  2.3} that (1.7) is equivalent to
$$a(z)v(w)=\sum_{j=0}^{N-1}(a_{(j)}v)(w){\dd_w^{j}}{{\delta(z-w)}\over j!},\eqno{(1.8)}$$
for some formal distributions $(a_{(j)}v)(w)$ with coefficients in $V$,
which
are uniquely determined by the formula
$$(a_{(j)}v)(w)={\rm Res}_{z}(z-w)^ja(z)v(w).$$
Formula (1.8) is obviously equivalent to
$$a_{[m]}v_{[n]}=\sum_{j\in\z_+}{m\choose j}(a_{(j)}v)_{[m+n-j]}.\eqno{(1.9)}$$

\noindent {\it Example 1.1.} Consider the following representation of the
Virasoro algebra $\V$ in the vector space $V$ with basis $v_{[n]}$,  
$n\in\Z$, over $\C$:
$$L_mv_{[n]}=((\Delta-1)(m+1)-n)v_{[m+n]}+\alpha v_{[m+n+1]},\eqno{(1.10)}$$
where $\alpha,\Delta\in\C$.  In terms of formal distributions $L(z)$ and
$v(z)$ this can be written as follows:
$$L(z)v(w)=(\dd_{w}+\alpha)v(w)\delta(z-w)+\Delta v(w)\dd_{w}\delta(z-w).\eqno{(1.11)}$$
Hence $L(z)$ and $v(z)$ are local.  This $\V$-module is denoted by $M^c(\alpha,\Delta)$ for reasons explained further.  A more invariant description of the module $M^c(\alpha,\Delta)$ is as follows.  Of course, $\V$ is the Lie algebra of regular vector fields on $\C^{\times}$, where $L_n=-t^{n+1}\dt$, $n\in\Z$.  For $\alpha,\lambda\in\C$ let
$$F_{\alpha,\lambda}=\C[t,t^{-1}]e^{-\alpha t}dt^{-\lambda}.$$
The Lie algebra $\V$ acts on the space $F_{\alpha,\lambda}$ in a natural way:
$$(g(t)\dt) f(t)dt^{-\lambda}=(g(t)f'(t)-\lambda f(t)g'(t))dt^{-\lambda}, \eqno{(1.12)}$$
where $g(t)\in\C[t,t^{-1}]$ and $f(t)\in \C[t,t^{-1}]e^{-\alpha t}$.  Then letting $v_{[n]}=t^ne^{-\alpha t}dt^{-\lambda}$ identifies $M^c(\alpha,\Delta)$ with $F_{\alpha,\Delta-1}$ (cf.~(1.10) and (1.12)).

Suppose that $V$ is spanned over $\C$ by the coefficients of a family $\E$ of  
formal distributions such that all $a(z)\in\F$ are local with respect to all  
$v(z)\in \E$.  Then we call $(V,\E)$ a {\it conformal module over $(\G,\F)$}.

One can show ([5]) that the family $\E$ of a conformal module $(V,\E)$ over
$(\G,\F)$ can always be included in a larger family $\bar{\E}$ which is still
local with respect to $\bar{\F}$, and such that $\dd \bar{\E}\subset \bar{\E}$ and
$a_{(j)}\bar{\E}\subset\bar{\E}$ for all $a\in\bar{\F}$ and $j\in\Z_+$.
It is straightforward to check the following properties for $a,b\in\bar{\F}$ and
$v\in\bar{\E}$ ($m,n\in\Z_+$):
$$[a_{(m)},b_{(n)}]v=\sum_{j=0}^m {m\choose j}(a_{(j)}b)_{(m+n-j)}v,\quad(\dd a)_{(n)} v=[\dd,a_{(n)}]v=-na_{(n-1)}v.\eqno{(1.13)}$$
(Here $[\cdot,\cdot]$ is the bracket of operators on $\bar{\E}$.)

Therefore it follows that any conformal module $(V,\E)$ over a Lie algebra of formal
distributions $(\G,\F)$ gives rise to a module $M=\bar{\E}$ over the conformal
algebra $R=\bar{\F}$, defined as follows.  It is a (left)
$\C[\dd]$-module equipped with a family of $\C$-linear maps $a\rightarrow
a_{(n)}^M$ of $R$ to ${\rm End}_{\c}M$, for each $n\in\Z_+$, such that the
following properties hold (cf.~(1.13)) for $a,b\in R$ and
$m,n\in\Z_+$:
\item{(M0)} $\quad a^M_{(n)}v=0$, for $v\in M$ and $n>>0$,
\item{(M1)} $\quad [a^M_{(m)},b^M_{(n)}]=\sum_{j=0}^m{m\choose
j}(a_{(j)}b)^M_{(m+n-j)}$,
\item{(M2)} $\quad (\dd a)^M_{(n)}=[\dd,a^M_{(n)}]=-na_{(n-1)}^M$.

Conversely, suppose that a conformal algebra $R=\oplus_{i\in
I}\C[\dd]a^i$ is a free
$\C[\dd]$-module and consider the associated Lie algebra of formal
distributions $(\G(R),\F)$.  Let $M$ be a module over the
conformal algebra $R$ and suppose that $M$ is a free  
$\C[\dd]$-module
with $\C[\dd]$-basis $\{v^{{\alpha}}\}_{{\alpha}\in J}$. This gives rise to
a conformal  module $M^c$ over $\G(R)$ with basis $v_{[n]}^{{\alpha}}$, where
$i\in I$, ${\alpha}\in J$ and $n\in\Z$, defined by (cf.~(1.8)):
$$a^{i}(z)v^{{\alpha}}(w)=\sum_{j\in\z_+}(a^{i}_{(j)}v^{{\alpha}})(w){\dd_w^{j}}{{\delta(z-w)}\over j!}.\eqno{(1.14)}$$
In the general case one defines $M^c$ as the quotient of the space $\oplus_{j\in\z}M_{[j]}$, where $M_{[j]}$ is a copy of $M$, by the $\C$-span of $\{(\dd v)_{[n]}+n v_{[n-1]}| v\in M, n\in\Z\}$.

Introducing the generating series $$a_{\lambda}b=\sum_{n=0}^{\infty}a_{(n)}b{{\lambda^n}\over {n!}}\quad {\rm and}\quad a_{\lambda}^M=\sum_{n=0}^{\infty}a_{(n)}^M{{\lambda^n}\over {n!}},$$ which lie in $R\otimes_{\c}\C[\lambda]$ and ${\rm End}_{\c}(M)\otimes_{\c}\C[\lambda]$ due to (C0) and (M0), respectively,  identities (M1) and (M2) can be written as
$$[a_{\lambda}^M,b_{\mu}^M]=(a_{\lambda}b)^M_{\lambda+\mu},\quad (\dd a)_{\lambda}^M=[\dd,a_{\lambda}^M]=-\lambda a_{\lambda}^M.\eqno{(1.15)}$$

\noindent{\it Example 1.2.} Let $M$ be a $\C[\dd]$-module. The following follows immediately from (1.15).
\item{(a)} A $R(\tilde{\G})$-module structure on $M$ is given by a $\C$-linear map $a\rightarrow a_{\lambda}^M$ of $\G$ to ${\rm End}_{\c}(M)[\lambda]$ such that
$$[a^M_{\lambda},b^M_{\mu}]=[a,b]_{\lambda+\mu}^M.\eqno{(1.16)}$$
\item{(b)} A $R(\V)$-module structure on $M$ is given by $L_{\lambda}^M\in {\rm End}_{\c}(M)[\lambda]$ such that
$$[L_{\lambda}^M,L_{\mu}^M]=(\lambda-\mu)L_{\lambda+\mu}^M.\eqno{(1.17)}$$
\item{(c)} A $R(\V)\lsemi R(\tilde{\G})$-module structure is given by $L_{\lambda}^M\in {\rm End}_{\c}(M)[\lambda]$ and $a_{\lambda}^M\in{\rm End}_{\c}(M)[\lambda]$ as in (b) and (c), satisfying (in addition to (1.16) and (1.17))
$$[L_{\lambda}^M,a_{\mu}^M]=-\mu a_{\lambda+\mu}^M.\eqno{(1.18)}$$

A conformal module $(V,\E)$ (respectively module $M$) over a Lie algebra  
of formal distributions $(\G,\F)$ (respectively over a conformal algebra  
$R$) is called {\it finite}, if $\bar{\E}$ (respectively $M$) is a finitely
generated $\C[\dd]$-module.

A conformal module $(V,\E)$ over $(\G,\F)$ is called {\it irreducible} if
there is no non-trivial invariant subspace which contains all $v_{[n]}$,
$n\in\Z$, for some non-zero $v\in \bar{\E}$.  Clearly a conformal module is  
irreducible if
and only if the associated module $\bar{\E}$ over the conformal algebra
$\bar{\F}$ is irreducible (in the obvious sense).

The above discussions reduce the
classification of finite conformal modules over a Lie algebra of formal  
distributions $(\G,\F)$ to the classification of finite modules over the
corresponding conformal algebra.

\noindent {\it Example 1.3.} (a) Example 1.1 gives the following family of $R(\V)$-modules $M(\alpha,\Delta)=\C[\dd]\otimes_{\c}\C v_{\Delta}$:
$$L_{(0)}v_{\Delta}=(\dd+{\alpha})v_{\Delta},\quad L_{(1)}v_{\Delta}=\Delta v_{\Delta},\quad L_{(j)}v_{\Delta}=0,\quad{\rm for}\ j\ge 2.\eqno{(1.19)}$$
\item{(b)} Any $\G$-module $(\pi,U)$ gives rise to a $R(\tilde{\G})$-module $M(U)=\C[\dd]\otimes_{\c}U$ defined by:
$$a_{(0)}u=\pi(a)u,\quad a_{(j)}u=0,\quad {\rm for\ } j\ge 1,\ {\rm where\ }a\in\G,\ u\in U.\eqno{(1.20)}$$
\item{(c)} Any $\G$-module $(\pi,U)$ and ${\alpha},\Delta\in\C$ gives rise to a $R(\V\lsemi\tilde{\G})$-module $M({\alpha},\Delta,U)=\C[\dd]\otimes_{\c}U$ defined by ($i\ge 0$, $j\ge 1$):
$$a_{(i)}u=\delta_{0,i}\pi(a)u,\quad L_{(0)}u=(\dd+{\alpha})u,\quad L_{(j)}u=\delta_{1,j}\Delta u.\eqno{(1.21)}$$

\noindent Note that the corresponding conformal modules in the cases (b) and (c) are $U\otimes_{\c}\C[t,t^{-1}]$ with the obvious action of $\tilde{\G}$ and the action of the Virasoro algebra given by (1.10), where $v_{[n]}=v\otimes t^n$ and $v\in U$.  As before, we will denote these corresponding conformal modules by $M^c(U)$ and $M^c(\alpha,\Delta,U)$, respectively.

\proclaim{Theorem 1.1 {\rm [1]}.} Let $\G$ be a simple Lie algebra.  The following is a complete list of finite irreducible, non-1-dimensional (over $\C$) modules over the conformal algebras $R(\V)$, $R(\tilde{\G})$ and $R(\V\lsemi\tilde{\G})$:
\item{(a)} $M({\alpha},\Delta)$, where $\Delta\not=0$.
\item{(b)} $M(U)$, where $U$ is a non-trivial finite-dimensional irreducible $\G$-module.
\item{(c)} $M({\alpha},\Delta,U)$, where either $U$ is a non-trivial finite-dimensional irreducible $\G$-module or $\Delta\not=0$.

In what follows we will always denote by $\G$ a finite-dimensional simple Lie algebra.  In the present paper we classify all extensions between two modules listed in Theorem 1.1 and between a module listed in Theorem 1.1 and a $1$-dimensional module. In the end we also state, without proofs, the classification of extensions between finite irreducible modules over the conformal algebra corresponding to the semidirect sum of the current algebra, associated with the $1$-dimensional Lie algebra, and the Virasoro algebra.

\noindent{\it Remark 1.1.} Given a module $M$ over a conformal algebra $R$, we may change its structure as a $\C[\dd]$-module by replacing $\dd$ by $\dd+A$, where $A$ is an operator on $M$ commuting with all the actions of $R$.  In particular if $A=\alpha\in\C$ we obtain a module over $R$ where the action of $\dd$ is replaced by $\dd+\alpha$.  Let us denote this module by $M_{\alpha}$. This twist by $\alpha$ for the corresponding conformal modules of Example 1.3 (b) and (c) has the following interpretation.  If $M(U)^c$ is the $\G(R)$-module $U\otimes\C[t,t^{-1}]$ with formal distributions $\{u(z)=\sum_{n\in\z}(u\otimes t^n)z^{-n-1}|u\in U\}$, then $M(U)^c_{\alpha}$ is the $\G(R)$-module $U\otimes\C[t,t^{-1}]e^{-\alpha t}$ with formal distributions $\{u^{\alpha}(z)=\sum_{n\in\z}(u\otimes t^ne^{-\alpha t})z^{-n-1}|u\in U\}$.  Note that in the case of the current algebra we have $U\otimes\C[t,t^{-1}]\cong U\otimes\C[t,t^{-1}]e^{-\alpha t}$ as $\tilde{\G}$-modules.  For the corresponding conformal modules of Example 1.3 (a) the effect of this twist corresponds to replacing $\C[t,t^{-1}]dt^{1-\Delta}$ by $\C[t,t^{-1}]e^{-\alpha t}dt^{1-\Delta}$.

Let $V$ and $W$ be two modules over a conformal algebra (or a Lie algebra) $R$. An extension of $W$ by $V$ is an exact sequence of $R$-modules of the form
$$0\ \longrightarrow\ V\ {\buildrel i\over\longrightarrow }\ {F}\ {\buildrel
p\over\longrightarrow }\ W\ \longrightarrow\ 0.$$ Two extensions
$0\longrightarrow V{\buildrel i\over\longrightarrow}{F}{\buildrel
p\over\longrightarrow}W\longrightarrow 0$ and $0\longrightarrow V{\buildrel
i'\over\longrightarrow}{F'}{\buildrel p'\over\longrightarrow}W\longrightarrow 0$
are said to be {\it equivalent} if there exists a commutative diagram of
of the form
$$\matrix{ 0&\longrightarrow&V&{\buildrel i\over\longrightarrow }&{F}&{\buildrel
p\over\longrightarrow }&W&\longrightarrow&0\cr
&&\smapdown{1_V}&&\smapdown{\psi}&&\smapdown{1_W}&&\cr
0&\longrightarrow&V&{\buildrel i'\over\longrightarrow }&{F'}&{\buildrel
p'\over\longrightarrow }&W&\longrightarrow&0,\cr }$$ where $1_V:V\longrightarrow
V$ and $1_W:W\longrightarrow W$ are the respective identity maps and
$\psi:F\longrightarrow F'$ is a homomorphism of modules.

The direct sum of modules $V\oplus W$ obviously gives rise to an
extension.  Extensions equivalent to it are called trivial extensions.  
In general an extension can be thought of as the direct sum
of vector spaces $E=V\oplus W$, where $V$ is a submodule of $E$, while for $w$ in $W$ we have:
$$a\cdot w=aw+\phi_a(w),\quad a\in R,$$
where $\phi_a:W\rightarrow V$ is a linear map satisfying the cocycle condition: $\phi_{[a,b]}(w)=\phi_a(bw)+a\phi_b(w)-\phi_b(aw)-b\phi_a(w)$, $b\in R$. The set of these
cocycles form a vector space over $\C$.  Cocycles equivalent to the trivial extension are called coboundaries.  They form a subspace and the quotient space by it is denoted by ${\rm Ext}(W,V)$.

The space $F^{+}_{\alpha,\lambda}=\C[t]e^{-\alpha t}dt^{-\lambda}$ has a natural action of the Lie algebra $\V_+\subset\V$ of regular vector fields on $\C$.  In [3] extensions between the trivial $\V_+$-module and $F^+_{\lambda}:=F^+_{0,\lambda}$, and between $F^+_{\lambda}$ and $F^+_{\mu}$ were considered.  The classification was achieved in a round about way using Goncharova's theorem.  The modules over the Virasoro conformal algebra that we consider in this paper are related to $F^+_{\lambda}$ as follows.  In [1] it was shown that the modules $M(0,\Delta)$ over the Virasoro conformal algebra correspond precisely to finite conformal irreducible modules $M_+(\Delta)$ over the ``annihilation subalgebra'' $\V_+$ of $\V$, namely: $F^+_{-\Delta}\cong M_+(\Delta)^*$, where $^*$ stands for the restricted dual.  Thus the results of this paper also give the answer to the extension problem considered in [3].  Our results agree with the results given in there, except for a few minor discrepancies that we would like to point out.  The cocycles for the extensions of the pairs $(\mu=-5,\lambda=0)$ and $(\mu={-7\pm\sqrt{19}\over 2},\lambda={5\pm\sqrt{19}\over 2})$ given in [3] appear to contain typos.  Also in the cases when $\lambda=\mu$ and $(\mu=-1,\lambda=0)$ our answers disagree slightly with the ones given there.  It appears that in either case there should be another non-trivial extension.

Note also that as $\V_+$-modules we have $M_+(\Delta)\cong F_{\Delta-1}/F^+_{\Delta-1}\cong (F^+_{-\Delta})^*$.

\beginsection{2. Extensions involving $1$-dimensional module}

Recall that every finite irreducible $R(\V)$-module is $M({\alpha},{\Delta})=\C[\dd]\otimes_{\c}\C v_{\Delta}$ defined by (1.19)
where ${\alpha}\in\C$ and $\Delta\not=0$, or else it is a $1$-dimensional module (Theorem 1.1a).

We first consider extensions of $R(\V)$-modules of the form
$$0\longrightarrow\C c_{\beta}\longrightarrow E\longrightarrow
M({\alpha,\Delta})\longrightarrow 0,\eqno{(2.1)}$$
where $\C c_{\beta}$ denotes the $1$-dimensional $R(\V)$-module with $L_{(n)}c_{\beta}=0$, for all $n$, and $\dd c_{\beta}=\beta c_{\beta}$.
We have $E=\C c_{\beta}\oplus\C[\dd] v_{\Delta}$.  The following identities hold in $E$:
$$L_{(0)}v_{{\Delta}}=\dd v_{{\Delta}}+\alpha v_{\Delta}+f_{0}c_{\beta},\qquad
L_{(1)}v_{{\Delta}}={\Delta} v_{{\Delta}}+f_{1}c_{\beta},\qquad
L_{(i)}v_{{\Delta}}=f_{i}c_{\beta},\quad i\ge 2,$$
where $f_{n}\in\C$ and all but finitely many of them are $0$.  Let $f(\lambda)=\sum_{n\ge 0}f_{n}{\lambda^n\over n!}$.
Thus $L_{\lambda}v_{{\Delta}}=\dd v_{{\Delta}}+\alpha v_{{\Delta}}+{\Delta}\lambda v_{{\Delta}}+f(\lambda)c_{\beta}.$  Since $E$ is a representation of $R(\V)$, applying both sides of (1.17) to $v_{\Delta}$, we obtain the following condition on the polynomial $f(\lambda)$: 
$$(\alpha+\beta+\lambda+{\Delta}\mu)f(\lambda)-(\alpha+\beta+\mu+{\Delta}\lambda)f(\mu)=(\lambda-\mu)f(\lambda+\mu).\eqno{(2.2)}$$
Put $\mu=0$ in (2.2) we obtain
$$(\alpha+\beta) f(\lambda)=(\alpha+\beta+\Delta\lambda)f(0).$$
Thus when $\alpha+\beta\not=0$, $f(\lambda)$ is a scalar multiple of $(\alpha+\beta+\Delta\lambda)$.  However, this polynomial corresponds to the trivial extension as follows: Let $E=\C[\dd]v_{{\Delta}}\oplus\C c_{\beta}$ be the direct sum of $R(\V)$-modules. Put $v_{\Delta}'=v_{\Delta}-c_{\beta}$.  Then a short calculation shows that $L_{\lambda}v_{\Delta}'=(\dd +\alpha+\Delta\lambda)v_{\Delta}' +({\alpha+\beta}+\Delta\lambda)c_{\beta}$.  Thus for the extension problem we may assume from now on that $\alpha+\beta=0$.  So (2.2) reduces to
$$(\lambda+{\Delta}\mu)f(\lambda)-(\mu+{\Delta}\lambda)f(\mu)=(\lambda-\mu)f(\lambda+\mu).\eqno{(2.3)}$$

So solving the extension problem above is equivalent to solving the functional equation (2.3) on the polynomial $f(\lambda)$.
We have the following

\proclaim{Proposition 2.1.}  There are non-trivial extensions of $R(\V)$-modules of the form (2.1)
if and only if $\alpha+\beta=0$ and ${\Delta}=1$ or $2$.  In these cases there exist unique (up to a scalar) non-trivial extensions, i.e.~${\rm dim}_{\c}({\rm Ext}(M(\alpha,1),\C c_{-\alpha}))={\rm dim}_{\c}({\rm Ext}(M(\alpha,2),\C c_{-\alpha}))=1$.

\proof
Since the functional equation (2.3) is homogeneous, we may assume that $f(\lambda)=a_n\lambda^n$, $a_n\in\C$.  We compare the coefficients of $\lambda^n\mu$ and $\mu^n\lambda$ on both sides and observe that
\item{(i)} There are no solutions for $n>3$.
\item{(ii)} If $n=2,3$ we must have $n-1={\Delta}$.
\item{(iii)} If $n=1$, then ${\Delta}$ can be arbitrary.
\item{(iv)} If $n=0$, then ${\Delta}=0$.

\noindent  Now since ${\Delta}\not=0$, case (iv) cannot occur.  Also the polynomial $f(\lambda)=a_1\lambda$ in case (iii) corresponds to the trivial extension as follows:  Let $E=M(\alpha,{\Delta})\oplus\C c_{\beta}$ be the trivial extension.  Set $v_{{\Delta}}'=v_{\Delta}-({a_1\over{\Delta}})c_{\beta}$.  Then $E=\C[\dd]v_{{\Delta}}'+\C c_{\beta}$ and $L_\lambda v_{{\Delta}}'=\sum_{n\ge 0}L_{(n)}v_{{\Delta}}'{\lambda^n\over n!}$ gives rise to the polynomial $a_1\lambda$.  Thus only the polynomials corresponding $f(\lambda)=a_2\lambda^2$ and $f(\lambda)=a_3\lambda^3$ give non-trivial extensions with ${\Delta}=1,2$, respectively.  $\qed$

\noindent {\it Remark 2.1.} The extension of conformal modules for ${\Delta}=1,2$ can be realized as follows.  For ${\Delta}=1$ 
consider the natural action of $\V$ on the space of functions on $S^1$ of the form $\C[t,t^{-1}]e^{-\alpha t}$.  This module is isomorphic to $M^c(\alpha,1)$. Let $\C c_{-\alpha}$ be the trivial $\V$-module.  Then $E^c=\C[t,t^{-1}]e^{-\alpha t}+\C c_{-\alpha}$, with $\V$ acting trivially on $c_{-\alpha}$ while acting on $\C[t,t^{-1}]e^{-\alpha t}$ via the formula $(f(t)\in\C[t,t^{-1}]$, $g(t)\in\C[t,t^{-1}]e^{-\alpha t})$
$$(f(t)\dt) g(t)=f(t)g'(t)+{\rm Res}_{t}(f''(t)g(t))c_{-\alpha}.$$  
For ${\Delta}=2$ consider the space of vector fields on $S^1$ of the form $\C[t,t^{-1}]e^{-\alpha t}dt^{-1}$.  It is isomorphic to $M^c(\alpha,2)$.  Letting $c_{-\alpha}$ as before we have $E^c=\C[t,t^{-1}]e^{-\alpha t}dt^{-1}+\C c_{-\alpha}$, with $\V$ acting as
$$(f(t)\dt) g(t)dt^{-1}=(f(t)g'(t)-f'(t)g(t))dt^{-1}+{\rm Res}_{t}(f'''(t)g(t))c_{-\alpha}.$$
The case $\alpha=0$ precisely gives the commutation relations of
the Virasoro algebra with center $\C c_0$.

Next we consider extensions of $R(\V)$-modules of the form
$$0\longrightarrow M({\alpha,\Delta})\longrightarrow E\longrightarrow
\C c_{\beta}\longrightarrow 0.\eqno{(2.4)}$$
As a vector space $E=\C[\dd]v_{\Delta}\oplus\C c_{\beta}$. We have for
$n\ge 0$ (in $E$): $$L_{(n)}c_{\beta}=f_n(\dd)v_{{\Delta}},\quad
f_n(\dd)\in\C[\dd].$$ Let $f(\dd,\lambda)=\sum_{n\ge
0}f_n(\dd){\lambda^n\over n!}$.  Applying both sides of (1.17) to
$c_{\beta}$ gives the functional equation:
$$(\dd+\alpha+{\Delta}\lambda)f(\dd+\lambda,\mu)-(\dd+\alpha+{\Delta}\mu)f(\dd+\mu,\lambda)=(\lambda-\mu)f(\dd,\lambda+\mu).$$
Replacing $\dd$ by $\bar{\dd}=\dd+\alpha$ and letting
$\bar{f}(\bar{\dd},\lambda)=f(\bar{\dd}-\alpha,\lambda)$ we can
rewrite the above equation in a homogeneous form
$$(\bar{\dd}+\Delta\lambda)\bar{f}(\bar{\dd}+\lambda,\mu)-(\bar{\dd}+\Delta\mu)\bar{f}(\bar{\dd}+\mu,\lambda)=(\lambda-\mu)\bar{f}(\bar{\dd},\lambda+\mu),\eqno{(2.5)}$$
so that we may assume that $\bar{f}(\bar{\dd},\lambda)$ is homogeneous
in $\bar{\dd}$ and $\lambda$.

We put $\mu=0$ in (2.5) and conclude that $\bar{f}(\bar{\dd},\lambda)$
is a scalar multiple of
$(\bar{\dd}+{\Delta}\lambda)g(\bar{\dd}+\lambda),$ where
$g(\bar{\dd}+\lambda)$ is a polynomial in $\bar{\dd}+\lambda$, if
degree of $\bar{f}$ is positive, or else
$\bar{f}(\bar{\dd},\lambda)=a_0$, $a_0\in\C$.  Now let $$\bar{\dd}
c_{\beta}=(\alpha+\beta) c_{\beta}+ a(\bar{\dd})v_{\Delta},$$ where
$a(\bar{\dd})$ is a polynomial in $\bar{\dd}$.  We use
$[\bar{\dd},L_{\lambda}]c_{\beta}=-\lambda L_{\lambda}c_\beta$ and
conclude that
$$(\lambda+\bar{\dd}-\alpha-\beta)\bar{f}(\bar{\dd},\lambda)=
a(\bar{\dd}+\lambda)(\bar{\dd}+\Delta\lambda).
\eqno{(2.6)}$$
It follows that $a(\bar{\dd})$ is a scalar multiple of
$(\bar{\dd}-\alpha-\beta)g(\bar{\dd})$.  This however corresponds to
the trivial extension.  Namely, let $E=M(\alpha,\Delta)\oplus\C
c_{\beta}$ be the direct sum of $R(\V)$-modules.  Set
$c_{\beta}'=c_{\beta}+g(\bar{\dd})v_{\Delta}$.  Then $$\eqalignno{
\bar{\dd} c_{\beta}'=&(\alpha+\beta) c_{\beta}+\bar{\dd}
g(\bar{\dd})v_{\Delta}=(\alpha+\beta)(c_{\beta}'-g(\bar{\dd})v_{\Delta})+
\bar{\dd}
g(\bar{\dd})v_{\Delta}\cr
=&(\alpha+\beta)
c_{\beta}'+(\bar{\dd}-\alpha-\beta)g(\bar{\dd})v_{\Delta},\cr
L_{\lambda}c_{\beta}'=&L_{\lambda}(c_{\beta}+g(\bar{\dd})v_{\Delta})=L_{\lambda}g(\bar{\dd})v_{\Delta}=g(\bar{\dd}+\lambda)(\bar{\dd}+\Delta\lambda)v_{\Delta}.\cr}$$

Thus for the extension problem we may assume that
$f(\dd,\lambda)=a_0$, $a_0\in\C$.  In this case (2.6) gives
$\alpha+\beta=0$, $\Delta=1$ and $f(\dd,\lambda)=a(\dd)=a_0$.  This
extension is non-trivial.

\proclaim{Proposition 2.2.} There are non-trivial extensions of $R(\V)$-modules of the form (2.4)
if and only if $\alpha+\beta=0$ and ${\Delta}=1$.  In this case there
exists a unique (up to a scalar) non-trivial extension.

\noindent {\it Remark 2.2.} The extension of the corresponding conformal modules in Proposition 2.2 can be obtained as follows. Consider the space of $1$-forms on $S^1$ of the form $\C[t,t^{-1}]e^{-\alpha t}dt$, which is isomorphic to $M^c(\alpha,0)$.  The exact $1$-forms form a $\V$-submodule isomorphic to $M^c(\alpha,1)$.  It follows that $\C[t,t^{-1}]e^{-\alpha t}dt/d(\C[t,t^{-1}]e^{-\alpha t})$ is the trivial $1$-dimensional $\V$-module.

Summarizing we have

\proclaim{Theorem 2.1.} All non-trivial extensions between a $1$-dimensional (over $\C$) module and a finite irreducible conformal module $M^c(\alpha,\Delta)$ over the Virasoro algebra are described by
Remarks 2.1 and 2.2.

We will now discuss current conformal algebras.  Let $\G$ be a finite-dimensional simple Lie algebra.  Let $\tilde{\G}$ and $R(\tilde{\G})$ denote its current and the corresponding (to $\tilde{\G})$ conformal algebra, respectively.  Let $(\pi,U)$ be a finite-dimensional irreducible (possibly trivial) $\G$-module.  Let $M(U)$ be the module $\C[\dd]\otimes U$ over $R(\tilde{\G})$ defined by (1.20).  A $1$-dimensional (over $\C$) module over the conformal algebra $R(\tilde{\G})$ is $\C c_{\alpha}$ with actions:
$$\dd c_{\alpha}=\alpha c_{\alpha},\quad a_{(n)}c_{\alpha}=0,\quad a\in\G, n\in\Z_+,\alpha\in\C.$$ 

First we look at extensions of $R(\tilde{\G})$-modules of the form
$$0\longrightarrow\C c_{\alpha}\longrightarrow E\longrightarrow
M(U)\longrightarrow 0.\eqno{(2.7)}$$
We have for $a\in\G$ and $u\in U$:
$$a_{(0)}u=\pi(a)u+\varphi^0_a(u)c_{\alpha},\quad
a_{(i)}u=\varphi^i_a(u)c_{\alpha},$$
where $\varphi_a^i:U\rightarrow\C$ are $\C$-linear maps.  As before we write down the generating series $a_{\lambda}u=\sum_{n\ge 0}a_{(n)}u{\lambda^n\over n!}$.  Applying both sides of (1.16) to $u$ we obtain the functional equation
$$\sum_{n\ge 0}\varphi_a^n(\pi(b)u){\lambda^n\over n!}-\sum_{n\ge 0}\varphi_b^n(\pi(a)u){\mu^n\over n!}=\sum_{n\ge 0}\varphi^n_{[a,b]}(u){{(\lambda+\mu)^n}\over n!}.$$
We compare the coefficient of $\lambda\mu^{n-1}$, $n\ge 2$, on both sides and see that $\varphi^n_{[a,b]}(u)=0$.  Since $u$ is arbitrary and $[\G,\G]=\G$ it follows that $\varphi^n_a=0$, for $n\ge 2$.

The coefficients of $1,\lambda,\mu$ give
\item{(i)}$\varphi^0_{a}(\pi(b)u)-\varphi^0_b(\pi(a)u)=\varphi^0_{[a,b]}(u)$.
\item{(ii)} $\varphi^1_a(\pi(b)u)=\varphi^1_{[a,b]}(u)$.
\item{(iii)} $-\varphi_b^1(\pi(a)u)=\varphi^1_{[a,b]}(u)$.

(ii) is clearly equivalent to (iii) and means that $\varphi^1:\G\rightarrow U^*$ is a $\G$-module homomorphism.  (i) on the other hand is equivalent to saying that $\varphi^0:\G\rightarrow U^*$ is a $1$-cocycle.  But ${\rm H}^1(\G;U^*)=0$, hence $\varphi^0$ corresponds to the trivial extension.

\proclaim{Proposition 2.3.} There are extensions of $R(\tilde{\G})$-modules of the form (2.7) if and only if $U\cong\G$ is the adjoint representation.  In this case the extension is unique up to a scalar.

\noindent{\it Remark 2.3.}  The corresponding extension of conformal modules over $\tilde{\G}$ of Proposition 2.3 can be realized as follows. Let $(\cdot|\cdot)$ be a non-degenerate symmetric invariant bilinear form on $\G$.  Identifying $U$ with $\G$ and $\C c_{-\alpha}$ with a $1$-dimensional module, we have $E^c=\G\otimes\C[t,t^{-1}]e^{-\alpha t}+\C c_{-\alpha}$. Then for $a,b\in\G$, $f(t)\in\C[t,t^{-1}]$ and $g(t)\in\C[t,t^{-1}]e^{-\alpha t}$ we have
$$(a\otimes f(t))\cdot(b\otimes g(t))=[a,b]\otimes f(t)g(t)+(a|b){\rm Res}_{t}(f'(t)g(t))c_{-\alpha},$$
where $\alpha\in\C$.  In the case $\alpha=0$ the extension is of course the well-known affine Kac-Moody algebra associated to the finite dimensional Lie algebra $\G$.

Now we will give a proof of

\proclaim{Proposition 2.4.} There are no non-trivial extensions of $R(\tilde{\G})$-modules of the form
$$0\longrightarrow M(U)\longrightarrow E\longrightarrow
\C c_{\alpha}\longrightarrow 0.\eqno{(2.8)}$$

\proof We have 
$$a_{(j)}c_{\alpha}=\sum_{i\ge 0}\dd^i\varphi^{i,j}(a),$$
where $\varphi^{i,j}:\G\rightarrow U$ are linear maps for $j\ge 0$.  We write the generating series
$$a_{\lambda}c_{\alpha}=\sum_{j\ge 0}a_{(j)}c_{\alpha}{\lambda^j\over j!}=\sum_{i,j\ge 0}\dd^i\varphi^{i,j}(a){\lambda^j\over j!}.$$

\proclaim{Lemma 2.1.} We can assume that $\varphi^{i,0}(a)=0$, $\forall i\in\Z_{+}$ and $\forall a\in\G$.

\proof Consider the subspace of $M(U)$
generated by $a_{(0)}c_{\alpha}$, $\forall a\in\G$.  If $N$ is the biggest integer such that $\varphi^{N,0}(a)\not=0$ for some $a\in\G$, then $a_{(0)}c_{\alpha}\subseteq U\oplus\dd U\oplus\cdots\oplus\dd^N U$.  Clearly this space is $\G$-invariant, hence the $\G$-module generated by $a_{(0)}c_{\alpha}$ is a finite-dimensional $\G$-module, say $W$. Restrict (2.8) to $W$ and we get
$$0\longrightarrow W\longrightarrow W+\C c_{\alpha}\longrightarrow
\C c_{\alpha}\longrightarrow 0,$$
which is an extension of $\G$-modules.  This corresponds to a $1$-cocycle of $\G$ with values in $W$, which is trivial.  Thus there exists a vector $w\in W$ such that 
$$\sum_{i}\dd^i\varphi^{i,0}(a)=a w,\quad\forall a\in\G.$$
We define an extension of $R(\tilde{\G})$-modules of the form
(2.8) by:
$$a_{(0)}c_{\alpha}=aw,\quad a_{(n)}c_{\alpha}=a_{(n)}w,\quad \dd c_{\alpha}=\alpha c_{\alpha}-\alpha w+\dd w,\qquad a\in\G, w\in W.$$
But this corresponds to the trivial extension.  Namely, let $c_{\alpha}'=c_{\alpha}-w$.  Then $\dd c_{\alpha}'=\alpha c_{\alpha}'$, and $a_{(n)}c_{\alpha}'=0$ and hence this defines a trivial cocycle.
Hence subtracting this $1$-cocycle from our original $1$-cocyle, if necessary, we may assume that $a_{(0)}c_{\alpha}=0$. $\qed$

Thus we may assume that
$$a_{\lambda}c_{\alpha}=\sum_{i\ge 0,j\ge 1}\dd^i\varphi^{i,j}(a){\lambda^j\over j!}.$$
Applying both sides of (1.16) to $c_{\alpha}$ gives the functional equation
$$\sum_{i\ge 0,j\ge 1}(\dd+\lambda)^i\pi(a)\varphi^{i,j}(b){\mu^j\over j!}- \sum_{i\ge 0,j\ge 1}(\dd+\mu)^i\pi(b)\varphi^{i,j}(a){\lambda^j\over j!}=\sum_{i\ge 0,j\ge 1}\varphi^{i,j}([a,b])\dd^i{{(\lambda+\mu)^j}\over j!}.$$
Now if $(\pi,U)$ is the trivial $\G$-module, it follows that $\varphi^{i,j}=0$, and so the $E$ decomposes.  Thus we may assume that $U$ is non-trivial.

\proclaim{Lemma 2.2.} If $U$ is non-trivial, then $\dd c_{\alpha}=\alpha c_{\alpha}$.

\proof We have $0=[\dd,a_{(0)}]c_{\alpha}=\dd a_{(0)}c_{\alpha}-a_{(0)}\dd c_{\alpha}$.  Thus $0=-a_{(0)}\dd c_{\alpha}$ by Lemma 2.1.  Hence $\dd c_{\alpha}-\alpha c_{\alpha}$ is a $\G$-invariant.  Since $U$ is irreducible and non-trivial, our assertion follows.  $\qed$

Let $M\ge 1$ be an integer. We compare the coefficient of $\lambda^{M-k}\dd^k$, $k\ge 1$, and get
$$-\pi(b)\varphi^{k,M-k}(a)=\varphi^{k,M-k}([a,b]).$$
This is to say that $\varphi^{k,M-k}:\G\rightarrow U$ is a $\G$-module homomorphism.  Thus we obtain

\proclaim{Corollary 2.3.} If $U\not\cong\G$, then there are no extensions of $R(\tilde{\G})$-modules of the form (2.8).

Hence we may assume, after identifying $\G$ with $U$, that
$$\varphi^{i,j}(a)=c_{i,j}a,\quad c_{i,j}\in\C.$$
Now by (1.15) $[\dd,a_{\lambda}]c_{\alpha}=-\lambda a_{\lambda}c_{\alpha}$.  Hence by Lemma 2.2 we have $\dd a_{\lambda}c_{\alpha}-\alpha a_{\lambda}c_{\alpha}=-\lambda a_{\lambda}c_{\alpha}$.  Comparing the coefficients of $\dd^{i}\lambda^n$ on both sides we obtain
$$-\alpha c_{0,n}=n c_{0,n-1},\quad n\ge 1,$$
$$c_{i-n,n}-\alpha c_{i-1,n}=-nc_{i,n-1},\quad i\ge 1,n\ge 1.$$
But we know that $c_{i,0}=0$ for $i\ge 0$.  Hence if $\alpha=0$, $c_{i-1,n}=-nc_{i,n-1}$, which implies that $c_{i,n}=0$ for all $i,n$.  If $\alpha\not=0$, then $c_{i,n}={1\over\alpha}(c_{i-1,n}+nc_{i,n-1})$.  By induction on $i+n$, it also follows that $c_{i,n}=0$. $\qed$

Translating back to the language of conformal modules over Lie algebras of formal distributions we obtain

\proclaim{Theorem 2.2.} Non-trivial extensions between a $1$-dimensional (over $\C$) $\tilde{\G}$-module and a conformal  $\tilde{\G}$-module $M^c(U)$, where $U$ is a finite-dimensional irreducible $\G$-module,  exist if and only if $U\cong\G$.  They are all described by Remark 2.3.

Finally we consider extensions of the semidirect sum $R(\V)\lsemi R(\tilde{\G})$.  Recall that the finite irreducible modules are either $M(\alpha,\Delta,U)$ or else it is a $1$-dimensional (over $\C$) module $\C c_{\beta}$ with $L_{(n)}c_{\beta}=a_{(n)}c_{\beta}=0$, for all $n\in\Z_+$, and $\dd c_{\beta}=\beta c_{\beta}$, where $\beta\in\C$.  Here $(\pi,U)$ is a finite-dimensional non-trivial irreducible $\G$-module or else if $U$ is trivial, then $\Delta\not=0$.  The conformal algebra acts on $M(\alpha,\Delta,U)=\C[\dd]\otimes U$ via ($a\in\G$ and $u\in U$):
$$a_{(0)}u=\pi(a)u,\quad L_{(0)}u=\dd u+\alpha u,\quad L_{(1)}u=\Delta u.$$

We first consider extensions of $R(\V\lsemi\tilde{\G})$-modules of the form
$$0\longrightarrow\C c_{\beta}\longrightarrow E\longrightarrow M(\alpha,\Delta,U)\longrightarrow 0.\eqno{(2.9)}$$
We may assume that $U$ is a non-trivial $\G$-module, for otherwise, due to Proposition 2.3, it reduces to extensions of the Virasoro conformal algebra already considered above. Restricting (2.9) to the action of $R(\tilde{\G})$ and using Theorem 2.2, we may assume that
$$\eqalignno{a_{\lambda}u=&\pi(a)u+\lambda\varphi_a(u)c_{\beta},\quad\qquad a_{\lambda}c_{\beta}=L_{\lambda}c_{\beta}=0,\cr
L_{\lambda}u=&\dd u+\alpha u+\Delta \lambda u+f_u(\lambda)c_{\beta},\quad \dd c_{\beta}=\beta c_{\beta},\cr}$$
where $\varphi_a(u)=0$ for $U\not\cong \G$, and $\varphi_a(u)=(a|u)$ is a symmetric invariant bilinear form on $\G$ for $U\cong\G$, and $f_u(\lambda)$ is a polynomial in $\lambda$.  Applying both sides of (1.18) to $u$ we obtain
$$-\mu(\lambda+\mu)\varphi_a(u)=f_{\pi(a)u}(\lambda)-(\mu+\alpha+\beta+\Delta\lambda)\mu\varphi_a(u).\eqno{(2.10)}$$
If $U\not\cong\G$, then $\varphi=0$.  From (2.10) it follows that $f_u=0$, for all $u\in U$, and hence the module decomposes.  Now suppose that $U\cong\G$ and $\varphi_a(u)=(a|u)\not=0$. (2.10) again implies that $f_u=0$, for all $u\in U$, and $\Delta=1$ and $\alpha+\beta=0$.  Thus we have

\proclaim{Proposition 2.5.} Non-trivial extensions of $R(\V\lsemi\tilde{\G})$-modules of the form (2.9) exist if and only if $\alpha+\beta=0$, $\Delta=1$ and $U\cong\G$.  Identifying $U$ with $\G$ the unique (up to a scalar) extension is given by ($a,b\in\G$):
$$a_{\lambda}b=[a,b]+\lambda(a|b)c_{-\alpha},\quad a_{\lambda}c_{-\alpha}=L_{\lambda}c_{-\alpha}=0,\quad L_{\lambda}a=\dd a+\alpha a+\lambda a,\quad \dd c_{-\alpha}=-\alpha c_{-\alpha}.$$

Next consider extensions of $R(\V\lsemi\tilde{\G})$-modules of the form
$$0\longrightarrow M(\alpha,\Delta,U)\longrightarrow E\longrightarrow\C c_{\beta}\longrightarrow 0.\eqno{(2.11)}$$
We can assume that $(\pi,U)$ is non-trivial, due to Proposition 2.4.  Using Theorem 2.2 we have ($a\in\G$, $u\in U$)
$$a_{\lambda}c_{\beta}=0,\quad a_{\lambda}u=\pi(a)u,\quad L_{\lambda}c_{\beta}=\sum_{i=1}^n f^i(\dd,\lambda)u_i,\quad L_{\lambda}u=(\dd+\alpha+\Delta\lambda)u,$$
where $\{u_1,u_2,\cdots u_n\}$ is $\C$-basis of $U$ and $f^i(\dd,\lambda)$ are polynomials in $\dd$ and $\lambda$. Applying both sides of (1.17) to $c_{\beta}$ we obtain
$$\sum_{i=1}^n(f^i(\dd+\lambda,\mu)(\dd+\alpha+\Delta\lambda)-
f^i(\dd+\mu,\lambda)(\dd+\alpha+\Delta\mu))u_i=\sum_{i=1}^n(\lambda-\mu)f^i(\dd,\lambda+\mu)u_i,$$
from which it follows that for $i=1,\cdots,n$
$$f^i(\dd+\lambda,\mu)(\dd+\alpha+\Delta\lambda)-f^i(\dd+\mu,\lambda) (\dd+\alpha+\Delta\mu)
=(\lambda-\mu)f^i(\dd,\lambda+\mu),$$
which is precisely (2.5).  Non-trivial solution to (2.5) exists only when $\Delta=1$ and  $f^i=a^i$, where $a^i\in\C$.  Now applying both sides of (1.18) to $c_{\beta}$ we obtain
$$\sum_{i=1}^n f^i(\dd+\mu,\lambda)\pi(a)u_i=0,\quad\forall a\in\G.$$
Thus $\pi(a)(\sum_{i=1}^n a^iu_i)=0$, for all $a\in\G$.  Since $U$ is non-trivial and $u_1,\cdots,u_n$ is a basis of $U$, it follows that $a^i=0$ for all $i$.  Thus the module decomposes.  Translating to conformal modules we have

\proclaim{Theorem 2.3.} Non-trivial extensions between a $1$-dimensional (over $\C$) $\V\lsemi \tilde{\G}$-module $\C c_{\beta}$ and a conformal module $M^c(\alpha,\Delta,U)$, where $U$ is a finite-dimensional irreducible non-trivial $\G$-module, exist if and only if $\alpha=\beta$, $U\cong\G$ and $\Delta=1$.  These are precisely the $\tilde{\G}$-modules of Remark 2.3, with completely reducible action of $\V$.

\beginsection{3. Extensions of Virasoro conformal modules}

In this section we will study extensions of $R(\V)$-modules of the form
$$0\longrightarrow M(\alpha,{\bar{\Delta}})\longrightarrow E\longrightarrow
M(\beta,{\Delta})\longrightarrow 0.\eqno{(3.1)}$$
As a $\C[\dd]$-module, we identify $E$ with the direct sum of $M(\alpha,\bar{\Delta})$ and $M(\beta,\Delta)$.
We have $L_{\lambda}v_{{\bar{\Delta}}}=\dd v_{{\bar{\Delta}}}+\alpha v_{{\bar{\Delta}}}+{{\bar{\Delta}}}\lambda v_{{\bar{\Delta}}}$ and $L_{\lambda}v_{{\Delta}}=\dd v_{{\Delta}}+\beta v_{{\Delta}}+{\Delta}\lambda v_{{\Delta}}+\sum_{n\ge 0}f_{n}(\dd){{\lambda^n}\over {n!}}v_{{\bar{\Delta}}}$, where $f_n(\dd)$ is a polynomial in $\dd$ for each $n$ and $f_n(\dd)=0$ for $n>>0$.  Putting $f(\dd,\lambda)=\sum_{n\ge 0}f_n(\dd){{\lambda^n}\over {n!}}$, we obtain, applying both sides of (1.17) to $v_{\Delta}$, the functional equation
$$\eqalignno{(\lambda-\mu)f(\dd,\lambda+\mu)=
&(\dd+\lambda+{\Delta}\mu+\beta)f(\dd,\lambda)+(\dd+\alpha+ {\bar{\Delta}}\lambda)f(\dd+\lambda,\mu)\cr
-&(\dd+\mu+{\Delta}\lambda+\beta)f(\dd,\mu)-(\dd+\alpha+ {\bar{\Delta}}\mu)f(\dd+\mu,\lambda).&{(3.2)}\cr}$$

Putting $\mu=0$ in (3.2) we get, after simplification,
$$(\beta-\alpha)f(\dd,\lambda)= (\dd+\Delta\lambda+\beta)f(\dd,0)-(\dd+\alpha+\bar{\Delta}\lambda)f(\dd+\lambda,0).$$
So if $\beta-\alpha\not=0$, then $f(\dd,\lambda)$ is a scalar multiple of $(\dd+\Delta\lambda+\beta)g(\dd)-(\dd+\alpha+\bar{\Delta}\lambda)g(\dd+\lambda)$, for some polynomial $g$.  But such a polynomial gives rise to the trivial extension as follows: Let $E=M(\alpha,\bar{\Delta})\oplus M(\beta,\Delta)$ and set $v_{\Delta}'=v_{\Delta}+g(\dd)v_{\bar{\Delta}}$.  Then a straightforward calculation shows that $L_{\lambda}v_{\Delta}'$ precisely gives rise to this polynomial.  Thus we may assume from now on that $\alpha=\beta$.  Next we put $\bar{\dd}=\dd+\alpha$ in (3.2) and let $\bar{f}(\bar{\dd},\lambda)=f(\bar{\dd}-\alpha,\lambda)$. Then the functional equation (3.2) becomes
$$\eqalign{(\lambda-\mu)\bar{f}(\bar{\dd},\lambda+\mu)=
&(\bar{\dd}+\lambda+{\Delta}\mu)\bar{f}(\bar{\dd},\lambda)+(\bar{\dd}+ {\bar{\Delta}}\lambda)\bar{f}(\bar{\dd}+\lambda,\mu)\cr
-&(\bar{\dd}+\mu+{\Delta}\lambda)\bar{f}(\bar{\dd},\mu)-(\bar{\dd} +{\bar{\Delta}}\mu)\bar{f}(\bar{\dd}+\mu,\lambda).\cr}$$
However, we will continue to write $\dd$ for $\bar{\dd}$ and $f$ for $\bar{f}$, but keep in mind that in the case when $\alpha\not=0$, we need to perform a shift by $\alpha$ in order to get the correct solutions in this case.  Thus we are to solve the functional equation
$$\eqalign{(\lambda-\mu)f(\dd,\lambda+\mu)=
&(\dd+\lambda+{\Delta}\mu)f(\dd,\lambda)+(\dd+ {\bar{\Delta}}\lambda)f(\dd+\lambda,\mu)\cr
-&(\dd+\mu+{\Delta}\lambda)f(\dd,\mu)-(\dd+ {\bar{\Delta}}\mu)f(\dd+\mu,\lambda).\cr}\eqno{(3.3)}$$

By the nature of (3.3) we may assume that a solution is a homogeneous polynomial in $\dd$ and $\lambda$ of degree $n$.  Hence we will assume from now on that $f(\dd,\lambda)=\sum_{j=0}^n a_j \dd^{n-j}\lambda ^j$, $a_j\in\C$.  The non-negative integer $n$ is called the {\it degree} of the extension.

\noindent {\it Remark 3.1.} The polynomial
$$(\dd+\lambda)^{n-1}(\dd+{\bar{\Delta}}\lambda)- (\dd+{{\Delta}}\lambda)\dd^{n-1}\eqno{(3.4)}$$
is always a solution to (3.3).  However this solution corresponds to the trivial extension for $n\ge 1$ as follows: Let $E=\C[\dd]v_{\bar{\Delta}}\oplus \C[\dd]v_{\Delta}$ and $v_{{\Delta}}, v_{{\bar{\Delta}}}$ be the corresponding generators.  For $n\ge 1$ we set $v_{{\Delta}}'=v_{{\Delta}}+\dd^n v_{{\bar{\Delta}}}$ and $v_{{\bar{\Delta}}}'=v_{{\bar{\Delta}}}$.  Then the polynomial $f(\dd,\lambda)$ corresponding to $E=\C[\dd]v_{{\Delta}}'\oplus\C[\dd]v_{{\bar{\Delta}}}'$ in the basis $v_{{\Delta}}'$ and $v_{{\bar{\Delta}}}'$ is precisely (3.4).

\proclaim{Lemma 3.1.} Let $f(\dd,\lambda)=\sum_{j=0}^n a_j \dd^{n-j}\lambda ^j$ be a solution to (3.3).  If $a_0\not=0$, then either
\item{(i)} $n=0$ and ${\bar{\Delta}}={\Delta}$, or else
\item{(ii)} $n=1$ with ${\bar{\Delta}}=0$ and ${\Delta}=1$.

\proof
Setting $\mu=0$ in (3.3) we obtain
$(\dd+{\bar{\Delta}}\lambda)f(\dd+\lambda,0)=(\dd+{\Delta}\lambda) f(\dd,0)$.
Hence
$a_0{\bar{\Delta}} (\dd+\lambda)^{n+1}+a_0({\Delta}-1)\dd^{n+1}=a_0({\bar{\Delta}}-1)\dd(\dd+\lambda)^n+a_0{\Delta} \dd^n(\dd+\lambda)$.  This gives, since $a_0\not=0$,
$${\bar{\Delta}} (\dd+\lambda)^{n+1}+({\Delta}-1)\dd^{n+1}=({\bar{\Delta}}-1)\dd(\dd+\lambda)^n+{\Delta} \dd^n(\dd+\lambda),$$
from which the assertion follows immediately. $\qed$

Solutions of small degrees can be checked directly.  We have

\proclaim{Proposition 3.1.} All solutions of (3.3) of degree less than or equal to 2 are as follows:
\item{(i)} $f(\dd,\lambda)=a_0$ for ${\Delta}={\bar{\Delta}}$.  It is non-trivial for $a_0\not=0$. 
\item{(ii)} $f(\dd,\lambda)=a_1\lambda$ for any ${\Delta},{\bar{\Delta}}\in\C$. It is equivalent to the trivial extension for ${\Delta}\not={\bar{\Delta}}$.  It is non-trivial in the case ${\Delta}={\bar{\Delta}}$ and $a_1\not=0$.
\item{(iii)} $f(\dd,\lambda)=a_0\dd+a_1\lambda$ for ${\Delta}=1$ and ${\bar{\Delta}}=0$.  It corresponds to a non-trivial extension if and only if $a_0\not=0$.
\item{(iv)} $f(\dd,\lambda)={{-a_2({\Delta}-{\bar{\Delta}}-1)}\over{{\bar{\Delta}}}}\dd\lambda+a_2\lambda^2$, for ${\Delta}$ arbitrary and ${\bar{\Delta}}\not=0$.  This is equivalent to the trivial extension.
\item{(v)} $f(\dd,\lambda)=a_1\dd\lambda$, for ${\Delta}\not=1$ and ${\bar{\Delta}}=0$.  This corresponds to the trivial extension.
\item{(vi)} $f(\dd,\lambda)=a_1\dd\lambda+a_2\lambda^2$, for ${\Delta}=1$ and ${\bar{\Delta}}=0$.  The corresponding extension is non-trivial unless $a_1=a_2=0$.

\proof The fact that these are all the solutions can be verified directly.  (i) is obvious.  Let $E=\C[\dd] v_{\Delta}\oplus \C[\dd] v_{\bar{\Delta}}$ be the trivial extension and let $v_{\Delta}$ and $v_{{\bar{\Delta}}}$ be the corresponding generating vectors.  For (ii) note that by the change of basis $v_{{\Delta}}'=v_{{\Delta}}+{{a_1}\over {({\bar{\Delta}}-{\Delta})}}v_{{\bar{\Delta}}}$ and $v_{{\bar{\Delta}}}'=v_{{\bar{\Delta}}}$ we obtain the polynomial solution of (ii).  Hence it is equivalent to the trivial extension when ${\Delta}\not={\bar{\Delta}}$.  (iii) then follows from (ii) and the fact that $f(\dd,\lambda)=a_0\dd$ obviously cannot be a trivial solution.  For (iv) we take $v_{{\Delta}}'=v_{{\Delta}}+{{2a_2}\over{{\bar{\Delta}}}}\dd v_{{\bar{\Delta}}}$ and $v_{{\bar{\Delta}}}'=v_{{\bar{\Delta}}}$ in order to see that it is indeed the trivial extension, while  for (v) we take $v_{{\Delta}}'=v_{{\Delta}}+{{a_1}\over {(1-{\Delta})}}\dd v_{{\bar{\Delta}}}$ and $v_{{\bar{\Delta}}}'=v_{{\bar{\Delta}}}$.  (vi) can be checked directly. $\qed$

\proclaim{Lemma 3.2.} Let $f(\dd,\lambda)=\sum_{j=0}^n a_j \dd^{n-j}\lambda ^j$ be a solution of (3.3).
If $n\ge 2$ and $n+{\bar{\Delta}}-{\Delta}-1\not=0$, then $f(\dd,\lambda)$ is a scalar multiple of (3.4).  If $n+{\bar{\Delta}}-{\Delta}-1=0$ and $n\ge 3$, then $a_1=0$.

\proof Differentiating (3.3) with respect to $\lambda$ we obtain
$$\eqalign{(\lambda&-\mu)f_{\lambda}(\dd,\lambda+\mu)+f(\dd,\lambda+\mu)=
(\dd+\lambda+{\Delta} \mu)f_{\lambda}(\dd,\lambda)+f(\dd,\lambda)\cr
+&(\dd+{\bar{\Delta}} \lambda)f_{\partial}(\dd+\lambda,\mu)+{\bar{\Delta}} f(\dd+\lambda,\mu)
-{\Delta} f(\dd,\mu) - (\dd+{\bar{\Delta}} \mu)f_{\lambda}(\dd+\mu,\lambda),\cr}\eqno{(3.5)}$$
where $f_{\partial}$ and $f_{\lambda}$ above denote the partial derivatives of $f(\dd,\lambda)$ with respect to $\dd$ and $\lambda$, respectively.
Now we put $\lambda=0$ and get
$$\eqalign{-\mu f_{\lambda}(\dd,\mu)+f(\dd,\mu)&=(\dd+{\Delta}\mu)f_{\lambda}(\dd,0)+f(\dd,0)+\dd f_{\partial}(\dd,\mu)+{\bar{\Delta}} f(\dd,\mu)\cr&-{\Delta} f(\dd,\mu)-(\dd+{\bar{\Delta}}\mu)f_{\lambda}(\dd+\mu,0).\cr}$$
Since $n\ge 2$ we may assume by Lemma 3.1 that $a_0=0$.  Also $f_{\lambda}(\dd,0)=a_1\dd^{n-1}$.  Combined with the fact that $\dd f_{\partial}(\dd,\lambda)+\lambda f_{\lambda}(\dd,\lambda)=n f(\dd,\lambda)$ we can simplify the above equation to
$$(n+{\bar{\Delta}}-{\Delta}-1)f(\dd,\mu)=a_1((\dd+\mu)^{n-1}(\dd+{\bar{\Delta}}\mu)- (\dd+{\Delta}\mu)\dd^{n-1}).$$
Now if $(n+{\bar{\Delta}}-{\Delta}-1)\not=0$, then $f(\dd,\lambda)$ is a scalar multiple of (3.4).  Since for $n\ge 3$ $(\dd+\mu)^{n-1}(\dd+{\bar{\Delta}}\mu)- (\dd+{\Delta}\mu)\dd^{n-1}$ is non-zero, the second assertion follows. $\qed$

Remark 3.1, Proposition 3.1 and Lemmas 3.1 and 3.2  allow us to assume from now on that $n\ge 3$ with $n+{\bar{\Delta}}-{\Delta}-1=0$, $a_0=a_1=0$.  Note that in this case we always have scalar multiples of (3.4) as non-zero solutions, which however correspond to the trivial extension.

\proclaim{Proposition 3.2.} If there exists a non-trivial extension of the form (3.1), then ${\Delta}-{\bar{\Delta}}$ is a non-negative integer.  Furthermore if ${\Delta}-{\bar{\Delta}}> 1$, then there exists at most one (up to a scalar) non-trivial extension.

\proof
The cases when $a_0\not=0$, $a_1\not=0$ or $n\le 2$ are dealt with in Proposition 3.1.
Hence we may assume that $a_0=a_1=0$ and $n-1={\Delta}-{\bar{\Delta}}$.
We set $\lambda=\mu$ in (3.5) and obtain
$$\eqalign{f(\partial,2\lambda)=&(\partial+({\Delta}+1)\lambda)f_{\lambda}(\partial,\lambda)+(1-{\Delta})f(\partial,\lambda)+{\bar{\Delta}} f(\partial+\lambda,\lambda)\cr
&+(\partial+{\bar{\Delta}} \lambda)f_{\partial}(\partial+\lambda,\lambda)-(\partial+{\bar{\Delta}} \lambda)f_{\lambda}(\partial+\lambda,\lambda).\cr}$$
We substitute $f(\dd,\lambda)=\sum_{j=0}^n a_j\dd^{n-j}\lambda$ into the equation above and, comparing coefficients on both sides of the resulting equation, we get 
$$\eqalign{(2^k-(k^2-k+&2))a_k=\sum_{j=0}^{k-1}(n-j){{n-j-1}\choose{k-j}}a_j+ {\bar{\Delta}}\sum_{j=0}^{k-1}(n-j){{n-j-1}\choose{k-j-1}}a_j\cr &-\sum_{j=1}^{k-1}j{{n-j}\choose{k-j+1}}a_j- {\bar{\Delta}}\sum_{j=1}^{k-1}(j-1){{n-j}\choose{k-j}}a_j,\quad 1\le k\le n-1,\cr}\eqno{(3.6)}$$
and
$$(2^n-(n^2-n+2))a_n={\bar{\Delta}}\sum_{j=0}^{n-1}(n-2j+1)a_j.\eqno{(3.7)}$$
Now $(2^k-(k^2-k+2))\not=0$ for $k\ge 4$.  Hence (3.6) and (3.7) are recursion formulas for $a_j$ with $j\ge 4$.  As $a_0=a_1=0$, it follows that every $a_j$, for $j\ge 4$ is determined by $a_2$ and $a_3$.  Since any scalar multiple of (3.4) is always a solution of (3.3), it follows that the space ${\rm Ext}(M(\alpha,\Delta),M(\alpha,\bar{\Delta}))$ is at most $1$-dimensional. $\qed$

We will continue to assume that $a_0=a_1=0$ and ${\Delta}-{\bar{\Delta}}=n-1$. Put $\lambda=-\mu$ in (3.3) and we obtain
$$(\dd+(1-{\Delta})\lambda)f(\dd,\lambda)+(\dd+{\bar{\Delta}} \lambda)f(\dd+\lambda,-\lambda)=(\dd+({\Delta}-1)\lambda)f(\dd,-\lambda)+(\dd-{\bar{\Delta}} \lambda)f(\dd-\lambda,\lambda),$$
from which it follows that $F(\dd,\lambda)=(\dd+(1-{\Delta})\lambda)f(\dd,\lambda)+(\dd+{\bar{\Delta}} \lambda)f(\dd+\lambda,-\lambda)$ is an even function with respect to $\lambda$.
Now the $\dd^{n+1-k}\lambda^k$-th coefficient of $F(\dd,\lambda)$ is
$$(1+(-1)^k)a_k+(1-{\Delta})a_{k-1}+\sum_{j=0}^{k-1}(-1)^j{{n-j}\choose{k-j}}a_j+{\bar{\Delta}}\sum_{j=0}^{k-1}(-1)^j{{n-j}\choose{k-j-1}}a_j.$$
For $F(\dd,\lambda)$ to be an even polynomial in $\lambda$ we must have for $1\le k\le n$ and $k$ odd
$$(1-{\Delta})a_{k-1}+\sum_{j=0}^{k-1}(-1)^j{{n-j}\choose{k-j}}a_j+{\bar{\Delta}}\sum_{j=0}^{k-1}(-1)^j{{n-j}\choose{k-j-1}}a_j=0.\eqno{(3.8)}$$
We will now use recursion formulas (3.6), (3.7) and (3.8) to prove the following

\proclaim{Proposition 3.3.} Non-trivial extensions of the form (3.1) with $\alpha=\beta=0$ of degree $n\ge 6$ occur if and only if $(\Delta,\bar{\Delta})$ is one of the following pairs:
$$(5,0),\qquad(1,-4),\qquad({{7\over 2}\pm{\sqrt{19}\over 2}},{-{5\over 2}\pm{\sqrt{19}\over 2}}).$$

\proof The main idea of the proof is to employ the above mentioned recursion formulas and derive a formula relating $a_2$ and $a_3$.  Now for the space of extensions to be $2$-dimensional, $a_2$, $a_3$ must be independent parameters.  Hence we obtain 2 equations involving ${\bar{\Delta}}$ and $n$ that must be satisfied simultaneously.  This gives the required restrictions.

By assumption ${\Delta}-{\bar{\Delta}}\ge 5$ and $n\ge 6$. From (3.6) and (3.7) we obtain the following formulas:
$$\eqalignno{a_4&={1\over 12}(n-2)(n-3)(n-4+3{\bar{\Delta}})a_2-{1\over 4}(n-3)(n-4+2{\bar{\Delta}})a_3,&(3.9a)\cr
a_5&={1\over {120}}(n-2)(n-3)(n-4)(n-5+4{\bar{\Delta}})a_2-{1\over 10}(n-4)(n-5+2{\bar{\Delta}})a_4,&(3.9b)\cr
32 a_6&={\bar{\Delta}}(3a_2+a_3-a_4-3a_5),&{\rm for\ } n=6,\quad(3.9c)\cr
32 a_6&={1\over 40}(n-2)(n-3)(n-4)(n-5)(n-6+5{\bar{\Delta}})a_2\cr
&+{1\over 24}(n-3)(n-4)(n-5)(n-6+4{\bar{\Delta}})a_3\cr
&-{1\over 6}(n-4)(n-5)(n-6+3{\bar{\Delta}})a_4- {3\over 2}(n-5)(n-6+2{\bar{\Delta}})a_5,&{\rm for\ }n\ge 7,\quad(3.9d)\cr
84a_7&={1\over 180}(n-2)(n-3)(n-4)(n-5)(n-6)(n-7+6{\bar{\Delta}})a_2\cr
&+{1\over 60}(n-3)(n-4)(n-5)(n-6)(n-7+5{\bar{\Delta}})a_3-{1\over 3}(n-5)(n-6)(n-7+3{\bar{\Delta}})a_5\cr
&-2(n-6)(n-7+2{\bar{\Delta}})a_6,&(3.9e)\cr
198a_8&={\bar{\Delta}}(5a_2+3a_3+a_4-a_5-3a_6-5a_7),&{\rm for\ } n=8,\quad(3.9f)\cr
198a_8&={1\over 1008}(n-2)(n-3)(n-4)(n-5)(n-6)(n-7)(n-8+7{\bar{\Delta}})a_2\cr
&+{1\over 240}(n-3)(n-4)(n-5)(n-6)(n-7)(n-8+6{\bar{\Delta}})a_3\cr
&+{1\over 120}(n-4)(n-5)(n-6)(n-7)(n-8+5{\bar{\Delta}})a_4\cr
&-{1\over 24}(n-5)(n-6)(n-7)(n-8+4{\bar{\Delta}})a_5-{1\over 2}(n-6)(n-7)(n-8+3{\bar{\Delta}})a_6\cr
&-{5\over 2}(n-7)(n-8+2{\bar{\Delta}})a_7,&{\rm for}\ n\ge 9.\quad(3.9g)\cr}$$

Now we look at (3.8) and for $k=7,9$ derive
$$\eqalignno{4a_6&={1\over 120}(n-2)(n-3)(n-4)(n-5)(n-6+5{\bar{\Delta}})a_2\hskip 1.9in\cr
&-{1\over 24}(n-3)(n-4)(n-5)(n-6+4{\bar{\Delta}})a_3\cr
&+{1\over 6}(n-4)(n-5)(n-6+3{\bar{\Delta}})a_4-{1\over 2}(n-5)(n-6+2{\bar{\Delta}})a_5,&(3.10a)\cr
6a_8&={1\over 5040}(n-2)(n-3)(n-4)(n-5)(n-6)(n-7)(n-8+7{\bar{\Delta}})a_2\cr
&-{1\over 720}(n-3)(n-4)(n-5)(n-6)(n-7)(n-8+6{\bar{\Delta}})a_3&(3.10b)\cr
&+{1\over 120}(n-4)(n-5)(n-6)(n-7)(n-8+5{\bar{\Delta}})a_4\cr
&-{1\over 24}(n-5)(n-6)(n-7)(n-8+4{\bar{\Delta}})a_5\cr
&+{1\over 6}(n-6)(n-7)(n-8+3{\bar{\Delta}})a_6 -{1\over 2}(n-7)(n-8+2{\bar{\Delta}})a_7.\cr}$$

First we consider the case $n=6$, i.e.~${\Delta}-{\bar{\Delta}}=5$.  We plug in (3.9a), (3.9b) and (3.9c) into (3.10a).  Using (3.9a) and (3.9b) to get rid of the terms $a_4$ and $a_5$ in the resulting equation, we obtain an equation involving only $a_2$ and $a_3$.  To be more precise we have the following equation
$$(24{\bar{\Delta}}+30{\bar{\Delta}}^2+6{\bar{\Delta}}^3)a_2-(24{\bar{\Delta}}+18{\bar{\Delta}}^2+3{\bar{\Delta}}^3)a_3=0.$$
Now in order to have a non-trivial extension $a_2$ and $a_3$ must be independent parameters.  Thus we must be able to solve the equations $(24{\bar{\Delta}}+30{\bar{\Delta}}^2+6{\bar{\Delta}}^3)=0$ and $(24{\bar{\Delta}}+18{\bar{\Delta}}^2+3{\bar{\Delta}}^3)=0$ simultaneously.  It can be checked that the only solutions are ${\bar{\Delta}}=0$ and ${\bar{\Delta}}=-4$.

For the case $n\ge 7$, i.e.~${\Delta}-{\bar{\Delta}}\ge 6$, we proceed as before.  Here we use (3.9a), (3.9b), (3.9d) and (3.10a).  We derive the following equations on $\bar{\Delta}$ and $n$:
$$\eqalign{(n-3)(n-4)&(n-5)(n^2+(4\bar{\Delta}-7)n+4\bar{\Delta}^2-8\bar{\Delta}+6)(n-2)(n-3+2\bar{\Delta})=0,\cr-(n-3)(n-4)&(n-5)(n^2+(4\bar{\Delta}-7)n+4\bar{\Delta}^2-8\bar{\Delta}+6)(n-2+3\bar{\Delta})=0.\cr}$$
Since $n\ge 7$, this gives $\bar{\Delta}={{(2-n)\pm\sqrt{3n-2}}\over 2}$.
In the case when $n=7$ we get ${\bar{\Delta}}={-{5\over 2}\pm{\sqrt{19}\over 2}}.$

The case of $n=8$ is checked via Equation (3.3).  The reason here is that the approach of $n=6,7$ above is not sufficient to rule out the existence of non-trivial extensions.  However, it can be verified directly, by plugging in a homogeneous polynomial into (3.3), that there are no solutions that is not a scalar multiple of (3.4).  To perform this calculation we have used the software package Maple V Release 3.

We now consider the case $n\ge 9$.  In this case we plug in (3.9a), (3.9b), (3.9d), (3.9e) and (3.9g) into (3.10b).  We again obtain an equation involving only $a_2$ and $a_3$, namely
$$(n-3)(n-4)(n-5)(n-6)(n-7)A(n,{\bar{\Delta}})((n-2)(n-3+3{\bar{\Delta}})a_2-3(n-2+2{\bar{\Delta}})a_3)=0,$$where
$$\eqalignno{A(n,{\bar{\Delta}})&=
48{\bar{\Delta}}^4+96(n-2){\bar{\Delta}}^3+(72n^2-386n+808){\bar{\Delta}}^2\cr
&+(24n^3-224n^2+968n-1232){\bar{\Delta}}+3n^4-44n^3+267n^2-946n+720.\cr}$$
As before, the coefficients of $a_2$ and $a_3$ must be $0$ in order for them to be independent parameters.  This gives two equations that must be simultaneously satisfied for $n\ge 9$:
$$\eqalignno{&(n-3)(n-4)(n-5)(n-6)(n-7)A(n,{\bar{\Delta}})(n-2)(n-3+3{\bar{\Delta}})=0,\cr
&(n-3)(n-4)(n-5)(n-6)(n-7)A(n,{\bar{\Delta}})(n-2+2{\bar{\Delta}})=0.\cr}$$
Since $n\ge 9$, we have $A(n,\bar{\Delta})=0$.  In addition the case of $n\ge 7$ considered above also tells us that $\bar{\Delta}={{(2-n)\pm\sqrt{3n-2}}\over 2}$.  So we are looking for integers $n\ge 9$ satisfying the equation $A(n,{{(2-n)\pm\sqrt{3n-2}}\over 2})=0$.  We have verified using Maple V Release 3 that the only integer solution is $n=1$.  $\qed$

We summarize the results of this section in the following two
theorems:

\proclaim{Theorem 3.1.} The complete list of solutions to (3.3), corresponding to non-trivial extensions of $R(\V)$-modules of the form (3.1) with $\alpha=\beta=0$, is given as follows, where in cases (iii)-(vii) $a_2,a_3\in\C$ are such that $f(\dd,\lambda)$ is not a multiple of (3.4):
\item{(i)} ${\Delta}={\bar{\Delta}}$ with ${\Delta}\in\C$.  $f(\dd,\lambda)=a_0+a_1\lambda$, $(a_0,a_1)\not=(0,0)$.
\item{(ii)} ${\Delta}=1$ and ${\bar{\Delta}}=0$.  $f(\dd,\lambda)=a_0\dd+b_0\dd\lambda+b_1\lambda^2$, where $(a_0,b_0,b_1)\not=(0,0,0)$.
\item{(iii)} ${\Delta}-{\bar{\Delta}}=2$, with ${\Delta}\in\C$.  $f(\dd,\lambda)=\lambda^2(a_2\dd+a_3\lambda)$.
\item{(iv)} ${\Delta}-{\bar{\Delta}}=3$, ${\Delta}\in\C$. $f(\dd,\lambda)=\lambda^2(a_2\dd^2+a_3\dd\lambda+{{\bar{\Delta}}\over 2}(a_2-a_3)\lambda^2)$.
\item{(v)} ${\Delta}-{\bar{\Delta}}=4$, ${\Delta}\in\C$. $f(\dd,\lambda)=\lambda^2(a_2\dd^3+a_3\dd^2\lambda+{1\over 2}((3{\bar{\Delta}}+1)a_2-(2{\bar{\Delta}}+1)a_3)\dd\lambda^2+{{\bar{\Delta}}\over 10}((1-3{\bar{\Delta}})a_2+(2{\bar{\Delta}}+1)a_3)\lambda^3)$.
\item{(vi)} ${\Delta}=5$ and ${\bar{\Delta}}=0$.  $f(\dd,\lambda)=a_2\lambda^2\dd^4+a_3\dd^3\lambda^3+(2a_2-{3\over 2}a_3)\dd^2\lambda^4+{{3a_3-2a_2}\over {10}}\dd\lambda^5$.
\item{(vi')} ${\Delta}=1$ and ${\bar{\Delta}}=-4$.
$f(\dd,\lambda)=a_2\lambda^2\dd^4+a_3\dd^3\lambda^3+({{9a_3-20a_2}\over 2})\dd^2\lambda^4+({{63a_3}\over 10}-17a_2)\dd\lambda^5+{2\over 5}(7a_3-20a_2)\lambda^6$.
\item{(vii)} ${\Delta}={{7\over 2}\pm{\sqrt{19}\over 2}}$ and ${\bar{\Delta}}={-{5\over 2}\pm{\sqrt{19}\over 2}}$. $f(\dd,\lambda)=({15\over 4}{\bar{\Delta}} a_2-{\bar{\Delta}} a_3+{33\over 28}a_2-{9\over 28}a_3)\lambda^7+(11{\bar{\Delta}} a_2+{7\over 2}a_2-3{\bar{\Delta}} a_3-a_3)\lambda^6\dd+(-3{\bar{\Delta}} a_3+11{\bar{\Delta}} a_2+{5\over 2}a_2)\lambda^5\dd^2+(-3a_3+5{\bar{\Delta}} a_2-2{\bar{\Delta}} a_3+5a_2)\lambda^4\dd^3+a_3\lambda^3\dd^4+a_2\lambda^2\dd^5$.

\proof (i), (ii) follows from Proposition 3.1. (iii)--(vii) follows from Proposition 3.2, Proposition 3.3 and the fact that these polynomials indeed satisfy Equation (3.3).  $\qed$

\proclaim{Theorem 3.2.} Non-trivial extensions of the form (3.1) exist only if $\alpha=\beta$.  Furthermore for each $\alpha\in\C$ the polynomials of these non-trivial extensions are obtained from the polynomials of Theorem 3.1 by replacing $\dd$ by $\dd+\alpha$.  The space ${\rm Ext}(M(\alpha,\Delta), M(\alpha, \bar{\Delta}))$ is
$2$-dimensional in case (i), $3$-dimensional in case (ii),
 $1$-dimensional in cases (iii)-(vii) of Theorem 3.1, and is trivial
for all other values of $\Delta$ and $\bar{\Delta}$.

In conclusion of this section we provide an explicit formula to translate the extensions in Theorem 3.2 into the language of conformal modules.  Recall that $M^c(\alpha,\bar{\Delta})\cong\C[t,t^{-1}]e^{-\alpha t}dt^{1-\bar{\Delta}}$. We will identify these modules.

Consider the extension of conformal modules
$$0\longrightarrow M^c(\alpha,\bar{\Delta})\longrightarrow E^c\longrightarrow M^c(\alpha,\Delta)\longrightarrow 0.$$
Then as a vector space $E^c\cong \C[t,t^{-1}]e^{-\alpha t}dt^{1-\bar{\Delta}}\oplus\C[t,t^{-1}]e^{-\alpha t}dt^{1-\Delta}$.  Now $\V$ acts on $\C[t,t^{-1}]e^{-\alpha t}dt^{1-\bar{\Delta}}$ as in (1.12), while on $\C[t,t^{-1}]e^{-\alpha t}dt^{1-\Delta}$ its action is given by (with $h(t)\in\C[t,t^{-1}]$ and $g(t)\in\C[t,t^{-1}]e^{-\alpha t}$):
$$\eqalignno{(h(t){\dt})(g(t)dt^{1-\Delta})&=
((1-\Delta)h'(t)g(t)+h(t)g'(t))dt^{1-\Delta}\cr
&-\sum_{i,k} a_{ik}(-1)^i (\sum_{j=0}^i{i\choose j}h^{(k+i-j)}(t)g^{(j)}(t))dt^{1-\bar{\Delta}},&(3.11)\cr}$$
where $f(\dd,\lambda)=\sum_{i,k}a_{ik}\dd^i\lambda^k$ is a polynomial in Theorem 3.1 determining the extension of the corresponding $R(\V)$-modules, and $h^{(i)}(t)$ denotes the $i$-th derivative with respect to $t$.  

To see this, we may assume that $\alpha=0$, due to Remark 1.1 and Theorem 3.2.  The first two summands on the right hand side of (3.11) are clear.  In order to see the last summand note that modulo $M^c(0,\Delta)$ we have
$$\eqalign{(-t^{m+1}\dt)(t^ndt^{1-\Delta})&=
\sum_{i,k}{{m+1}\choose k}k!a_{ik}(\dd^iv_{\bar{\Delta}})_{[m+n+1-k]}\cr
&=\sum_{i,k}{{m+1}\choose k}k!a_{ik}(-1)^i{{m+n+1-k}\choose{i}}i!t^{m+n+1-k-i} dt^{1-\bar{\Delta}}.\cr}$$
But now one can show that
$${{m+1}\choose k}{{m+n+1-k}\choose i}k!i!t^{m+n+1-k-i}=-\sum_{j=0}^i{i\choose j}(\dt)^{k+i-j}(-t^{m+1})(\dt)^{j}(t^n).$$
Using this identity (3.11) follows.

\beginsection{4. Extensions of current conformal modules}

Throughout this section we let $\G$ be a finite-dimensional simple Lie algebra.  In this section we consider extensions of $R(\tilde{\G})$-modules of the form
$$0\longrightarrow M(V)\longrightarrow E\longrightarrow
M(U)\longrightarrow 0.\eqno{(4.1)}$$

First we study the case when $M(U)=\C[\dd]U$ and $M(V)=\C[\dd]V$, where
$(\pi,U)$ and $({\rho},V)$ are non-trivial irreducible representations of $\G$.  We shall identify $E$ with the direct sum of $\C[\dd]$-modules $M(U)$ and $M(V)$. 
We have for $a\in\G$, $v\in V$ and $u\in U$:
$$a_{\lambda}v=\rho(a)v,\quad a_{\lambda}u=\pi(a)u+\sum_{j,k\ge 0}\lambda^{j}\dd^k\varphi^{j,k}_a
u,\eqno{(4.2)}$$
where $\varphi^{j,k}:\G\otimes
U\rightarrow V$ are $\C$-linear maps.
Let us write $a_{(0)}u=\pi(a)u+\Phi_a(u)$, where $\Phi_a:U\rightarrow\C[\dd]V$ is a linear map.  This gives rise to a $\C$-linear map $\Phi:\G\rightarrow U^*\otimes\C[\dd]V$.  Since $\Phi$ is the restriction of our original $1$-cocycle to $\G$, it is a $1$-cocycle of $\G$.  However ${\rm Im}\Phi$ is contained in a finite-dimensional $\G$-module and therefore, due to the fact that $H^{1}(\G;W)=0$ for any finite-dimensional $\G$-module $W$, there exists $\omega\in U^*\otimes\C[\dd]V$ such that $\Phi_a=a\cdot\omega$.  
This $\omega$, as in the proof of Lemma 2.1, allows us to define an extension, equivalent to the trivial extension, and
therefore, replacing $u$ by $u-\omega(u)$ if necessary, we may assume that $a_{(0)}u=\pi(a)u$. This reduces (4.2) to
$$a_{\lambda}u=\pi(a)u+\sum_{j\ge 1,k\ge 0}\lambda^{j}\dd^k\varphi^{j,k}_a
u.$$
Applying both sides of (1.16) to $u$ we obtain, the following identity:
$$\eqalign{\sum_{j,k}\lambda^j\dd^k\varphi^{j,k}_a\pi(b)u&+\sum_{j,k}\mu^j(\dd+\lambda)^k{\rho}(a)\varphi^{j,k}_b
u -\sum_{j,k}\mu^j\dd^k\varphi^{j,k}_b\pi(a)u\cr
&-\sum_{j,k}\lambda^j(\dd+\mu)^k{\rho}(b)\varphi^{j,k}_a u=
\sum_{j,k}(\lambda+\mu)^{j}\dd^k\varphi^{j,k}_{[a,b]}u,}\eqno{(4.3)}$$
where $j\ge 1$ and $k\ge 0$. We will now proceed to study identity (4.3) in detail.  Before doing so we need
the following

\proclaim{Lemma 4.1.} Let $(\pi,U)$ and
$({\rho},V)$ be two non-trivial irreducible representations of $\G$.  Let
$T:\G\otimes U\rightarrow V$ be a homomorphism of $\G$-modules.  Suppose that
${\rho}(a)T(b\otimes u)+{\rho}(b)T(a\otimes u)=0$, $\forall a,b\in\G$.  Then $T=0$.

\proof Let $\H$ be a Cartan subalgebra of $\G$. Let 
${\alpha}_1,\cdots,{\alpha}_l$ be a set of simple roots.  Let $e_i\in\G_{{\alpha}_i}$
and $f_i\in\G_{-{\alpha}_i}$.  Let $u_0$ and $v_0$ be highest weight vectors of
$U$ and $V$, respectively.

Define $T_a:U\rightarrow V$ by $T_{a}(u)=T(a\otimes u)$.  Since $T$ is a
$\G$-homomorphism, we have ${\rho}(a)T_b(u)=T_b\pi(a)(u)+T_{[a,b]}(u)$, for all
$a,b\in\G$ and $u\in U$.  This together with the hypothesis of the lemma gives
$T_a\pi(b)(u)+T_b\pi(a)(u)=0$ for all $a,b\in\G$. Thus in particular
$T_a\pi(e_i)u_0+T_{e_i}\pi(a)u_0=0$.  Hence $T_{e_i}\pi(a)u_0=0$ for all
$a\in\G$.  In particular, $T_{e_i}\pi(h)u_0=0$ for all $h\in\H$.  Since $U$ is
non-trivial, it follows that
$$T_{e_i}u_0=0,\quad \forall i=1,\cdots,l.\eqno{(4.4)}$$

Now ${\rho}(a)T_{e_i}u_0+{\rho}(e_i)T_{a}u_0=0, \forall a\in\G$.  So by (4.4) we have 
${\rho}(e_i)T_{a}u_0=0$, $\forall a\in\G$.

Since $T$ is a $\G$-module homomorphism, we have
$0={\rho}(e_i)T_{a}u_0=T_{[e_i,a]}u_0+T_a\pi(e_i)u_0$.  Thus $T_{[e_i,a]}u_0=0$,
$\forall a\in\G$.  Therefore, in particular:
$$T_h u_0=0,\quad\forall h\in \H.\eqno{(4.5)}$$

We will now show that $T_{a}u_0=0$, for all $a\in\G$, hence
$T(\G\otimes u_0)=0$.  From this it follows readily that
$T=0$, since $T$ is a $\G$-homomorphism.

We have
${\rho}(a)T_{e_i} u_0+{\rho}(e_i)T_{a} u_0=0$.  So by (4.4) we
have $T_a u_0=c v_0$, for some $c\in\C$.  Now
${\rho}(a)T_h u_0+{\rho}(h)T_{a} u_0=0$ together with (4.5)
gives $c{\rho}(h)v_0=0$, for all $h\in\H$.  But $V$ is non-trivial, hence $c=0$.
$\qed$

We now turn our attention to (4.3).  (4.3) is a polynomial identity in the
variables
$\dd,\lambda,\mu$.  We obtain the following identities by comparing the
coefficients of the polynomials of degree 1 and 2:
$$\varphi^{1,0}_a\pi(b)-{\rho}(b)\varphi^{1,0}_a=\varphi_{[a,b]}^{1,0},\quad\varphi^{2,0}_a\pi(b)-{\rho}(b)\varphi^{2,0}_a=\varphi_{[a,b]}^{2,0},\quad
{\rho}(a)\varphi^{1,1}_b-{\rho}(b)\varphi^{1,1}_a=2\varphi_{[a,b]}^{2,0}.$$
The first two identities mean that $\varphi^{1,0}$ and $\varphi^{2,0}$ are
$\G$-module homomorphisms.  We will come back to the third identity later.

Let us now compare the coefficients of polynomials of degree $M$ in (4.3), where
$M\ge 3$.  Comparing the coefficients of $\lambda^{M-k}\dd^k$,
$\lambda^{M-j}\mu^j$,
$\lambda^{M-j-1}\mu^j\dd$ and $\lambda^i\mu^i\dd^{M-2i}$, respectively we obtain
the following identities:
\item{(a)}$\varphi^{M-k,k}_a\pi(b)-{\rho}(b)\varphi_a^{M-k,k}=
\varphi^{M-k,k}_{[a,b]},\quad 0\le k\le M-1.$
\item{(b)}${\rho}(a)\varphi_b^{j,M-j}-{\rho}(b)\varphi_a^{M-j,j}={M\choose
j}\varphi^{M,0}_{[a,b]},\quad 1\le j\le M-1.$
\item{(c)}$(M-j){\rho}(a)\varphi_b^{j,M-j}-(j+1){\rho}(b)\varphi_a^{M-j-1,j+1}={{M-1}\choose
j}\varphi_{[a,b]}^{M-1,1},\quad 1\le j\le M-2.$
\item{(d)}
${{M-i}\choose i}{\rho}(a)\varphi_b^{i,M-i}-{{M-i}\choose i}
{\rho}(b)\varphi_a^{i,M-i}={{2i}\choose i}\varphi^{2i,M-2i}_{[a,b]},
\quad 1\le i < {M\over 2}.$

\proclaim{Proposition 4.1.} For $M\ge 3$ we have
\item{(i)} $\varphi^{M-k,k}:\G\otimes U\rightarrow V$ is a $\G$-homomorphism for
$0\le k\le M-1$.
\item{(ii)} $\varphi^{j,M-j}=\varphi^{M-j,j},\quad 1\le j\le M-1$.
\item{(iii)} $\varphi^{M-j,j}={M\choose j}\varphi^{M,0},\quad 1\le j\le M-1$.
\item{(iv)} ${\rho}(b)\varphi_a^{M,0}=\varphi^{M,0}_b\pi(a)$.

\proof (i) is precisely the statement of (a).  Now interchanging $a$ and $b$ in
(b) above we get
${\rho}(b)\varphi_a^{j,M-j}-{\rho}(a)\varphi_b^{M-j,j}=-{M\choose
j}\varphi^{M,0}_{[a,b]}$.  Adding (b) to it we get
${\rho}(a)(\varphi_b^{j,M-j}-\varphi_b^{M-j,j})+
{\rho}(b)(\varphi_a^{j,M-j}-\varphi_a^{M-j,j})=0.$ It follows from Lemma 4.1 that
$\varphi^{j,M-j}=\varphi^{M-j,j}$, proving (ii). Next we interchange $a$ and $b$
in (c) and get 
$(M-j){\rho}(b)\varphi_a^{j,M-j}-(j+1){\rho}(a)\varphi_b^{M-j-1,j+1}=-{{M-1}\choose
j}\varphi_{[a,b]}^{M-1,1}$.  As before we add (c) to this identity and arrive at
$${\rho}(a)((M-j)\varphi_b^{j,M-j}-(j+1)\varphi_b^{M-j-1,j+1})+
{\rho}(b)((M-j)\varphi_a^{j,M-j}-(j+1)\varphi_a^{M-j-1,j+1})=0.$$ Using Lemma 4.1
again we obtain
$(M-j)\varphi^{j,M-j}=(j+1)\varphi^{M-j-1,j+1}$.  Thus by (ii)
$(M-j)\varphi^{M-j,j}=(j+1)\varphi^{M-j-1,j+1}$, and hence
$\varphi^{M-j,j}={M\choose j}\varphi^{M,0}$, for $0\le j\le M-1$. Finally
substituting (iii) into (d) we obtain
${{M-i}\choose i}{M\choose i}({\rho}(a)\varphi_b^{M,0}-{\rho}(b)\varphi_a^{M,0})=
{2i\choose i}{M\choose{2i}}\varphi_{[a,b]}^{M,0}$, which is equivalent to
${\rho}(a)\varphi_b^{M,0}-{\rho}(b)\varphi_a^{M,0}=\varphi_{[a,b]}^{M,0}$. 
Replacing $\varphi_{[a,b]}^{M,0}$ by 
${\rho}(a)\varphi_b^{M,0}-\varphi_b^{M,0}\pi(a)$ (due to (a)), (iv) follows.
$\qed$

An immediate consequence of Proposition 4.1 (iii) is the following

\proclaim{Corollary 4.1.} For $M\ge 3, j\ge 1, k\ge 0$ we have
$$\sum_{j+k=M}\lambda^j\dd^k\varphi^{j,k}=
((\dd+\lambda)^{M}-\dd^M)\varphi^{M,0}.$$ In particular $\varphi^{j,M-j}$ is
determined by $\varphi^{M,0}$.

We will next study the identity (iv) in Proposition 4.1.

\proclaim{Lemma 4.2.} Let $(\pi,U)$ and
$({\rho},V)$ be non-trivial irreducible representations of $\G$.  Let $T:\G\otimes U\rightarrow V$ be a
$\G$-homomorphism such ${\rho}(a)T(b\otimes u)=T(a\otimes\pi(b)u)$, for all
$a,b\in\G$.
\item{(i)} If $(\pi,U)\not\cong ({\rho},V)$, then $T=0$.
\item{(ii)} If $(\pi,U)\cong ({\rho},V)$, then there exists $c_0\in\C$ such that
$T(a\otimes u)=c_0\pi(a)u$, for all $a\in\G,u\in U$.

\proof We will continue to use the notation in the proof of Lemma 4.1.  Let
furthermore $\Lambda$ and $\Lambda'$ denote the highest weights of the
representations $\pi$ and ${\rho}$, respectively.  For $a\in\G$ we have
${\rho}(a)T_{e_i}u_0=T_a\pi(e_i)u_0$ by hypothesis.  Since $u_0$ is a highest
weight vector, we have ${\rho}(a)T_{e_i}u_0=0$, for any $a\in\G$.  But $U$ is
non-trivial, hence $T_{e_i}u_0=0$.

Now for $h\in\H$, ${\rho}(e_i)T_{h}u_0=T_{e_i}\pi(h)u_0={\Lambda(h)}T_{e_i}u_0=0$. 
Thus $T_h u_0=c(h) v_0$, for some $c\in\H^*$.

Suppose that $U\not\cong V$. Since $T$ is a $\G$-homomorphism we also have
${\rho}(h')T_h u_0=T_h\pi(h')u_0$ for $h,h'\in\H$.  Hence
$\Lambda'(h')c(h)v_0=\Lambda(h')c(h)v_0$, which implies that
$(\Lambda(h')-\Lambda'(h'))c(h)=0$, $\forall h,h'\in\H$ and hence $c=0$. 
Therefore
$T_h u_0=0$. Thus $\Lambda(h)T_a u_0=T_a\pi(h)u_0={\rho}(a)T_h u_0=0$, $\forall
a\in\G,h\in\H$.  But $\Lambda\not=0$, hence $T_a u_0=0$, $\forall a\in\G$.  It
follows that $T=0$, proving (i).

Suppose that $U\cong V$.  By assumption $\pi(h)T_{h'}u_0=T_h\pi(h')u_0$, $\forall
h,h'\in\H$, which is to say that $\Lambda(h)c(h')u_0=\Lambda(h')c(h)u_0$,
$\forall h,h'\in\H$.  Since $\Lambda\not=0$, there exists $h_0\in\H$ such that
$\Lambda(h_0)\not=0$.  Thus $c(h)={{c(h_0)}\over {\Lambda(h_0)}}\Lambda(h)$, for
all $h\in\H$.  Set ${{c(h_0)}\over {\Lambda(h_0)}}=c_0$, and we have
$c(h)=c_0\Lambda(h)$ and thus
$T_hu_0=c_0\Lambda(h)u_0$, $\forall h\in\H$.  Thus
$\Lambda(h)T_a u_0=T_a\pi(h)u_0=\pi(a)T_h u_0=c_0\pi(a)\Lambda(h) u_0$, $\forall
a\in\G,h\in\H$.  Since $\Lambda\not=0$, this gives $T_a u_0=c_0\pi(a) u_0$,
$\forall a\in\G$. (ii) follows then from the fact that $T$ is a $\G$-module
homomorphism. $\qed$

From Lemma 4.2 and Proposition 4.1 (iv) we obtain immediately.

\proclaim{Corollary 4.2.} Let $M\ge 3$.
\item{(i)} If $(\pi,U)\not\cong(\rho,V)$, then $\varphi^{M,0}=0$.
\item{(ii)} If $(\pi,U)\cong(\rho,V)$, then $\varphi^{M,0}=c_M\pi(a)$, $\forall
a\in\G$ and $c_M\in\C$.

\proclaim{Proposition 4.2.}
\item{(i)} If $({\rho},V)\not\cong(\pi,U)$, then
$a_{\lambda}u=\pi(a)u+\lambda\varphi^{1,0}(a\otimes
u)+\lambda^2\varphi^{2,0}(a\otimes u)+2\lambda\dd\varphi^{1,1}(a\otimes u)$
\item{(ii)} If $({\rho},V)\cong(\pi,U)$, then
$a_{\lambda}u=\pi(a)u+\lambda\varphi^{1,0}(a\otimes
u)+\lambda^2\varphi^{2,0}(a\otimes u)+2\lambda\dd\varphi^{1,1}(a\otimes u)
+\sum_{M\ge 3}c_M((\dd+\lambda)^{M}-\dd^M)\pi(a)u$, for some $c_M\in\C$.
\item{(iii)} In (i) and (ii) above we have
${\rho}(a)\varphi^{1,1}_b-{\rho}(b)\varphi_a^{1,1}=\varphi^{2,0}_{[a,b]}$, $\forall
a,b\in\G$.

\noindent{\it Remark 4.1.} The formula in Proposition 4.2 (ii) above
can be brought to the form
$$a_{\lambda}u'=\pi(a)u'+\lambda\varphi^{1,0}(a\otimes u')
+\lambda^2\varphi^{2,0}(a\otimes u')+2\lambda\dd\varphi^{1,1}(a\otimes u'),$$ by making a suitable choice of $U'\subseteq\C[\dd](U\oplus V)$ such that $U'\cong U$ (as $\G$-modules) and $\C[\dd](U\oplus V)=\C[\dd](U'\oplus V)$.  Explicitly we set $u'=u-\sum_{M}c_M\dd^M \bar{u}$, where $\bar{u}$ here is to denote that it is an element of $V$ corresponding (via an isomorphism of $\G$-modules) to $u$, and leave $V$ unchanged.

\noindent {\it Remark 4.2.} Here we need to single out the $sl_2$ case.  In this
case we can view Proposition 4.2 (iii) above as the definition of $\varphi^{2,0}$.  It can be
readily checked (using the fact that $\varphi^{1,1}$ is a $\G$-homomorphism)
that, in the case of $sl_2$, $\varphi^{2,0}$ is indeed a $\G$-module homomorphism.

This discussion leads to the classification for $\G=sl_2$.

\proclaim{Theorem 4.1.} Non-trivial extensions of the form (4.1) for $\G=sl_2$ with $(\pi,U)$ and $({\rho},V)$
non-trivial irreducible $\G$-modules exist only if $(\pi,U)\not\cong({\rho},V)$.  They are given by:
$$a_{\lambda}v=\rho(a)v,\quad a_{\lambda}u=\pi(a)u+\lambda\varphi^{1,0}(a\otimes
u)+\lambda^2\varphi^{2,0}(a\otimes u)+2\lambda\dd\varphi^{1,1}(a\otimes u),$$
where $a\in\G$, $u\in U$, $v\in V$ and $\varphi^{1,0},\varphi^{1,1}:\G\otimes U\rightarrow V$ are arbitrary
homomorphisms of $\G$-modules, and $\varphi^{2,0}$ is determined by
$${\rho}(a)\varphi^{1,1}_b-{\rho}(b)\varphi_a^{1,1}=\varphi^{2,0}_{[a,b]},\quad \forall a,b\in\G.\eqno{(4.6)}$$

\proof If $U\not\cong V$, it follows from Proposition 4.2 (i) and Remark 4.2 that all non-trivial extensions are of form above.  Now if $U\cong V$ note
that the map
$\varphi:\G\otimes U\rightarrow U$ defined by $\varphi(a\otimes u):=\pi(a)u$ is a homomorphism of $\G$-modules.  Since $U$ appears with multiplicity 1
in the $\G\otimes U$ it follows that when $(\pi,U)\cong({\rho},V)$
$\varphi^{1,0}(a\otimes u)=c_1\pi(a)u$ and  $\varphi^{1,1}(a\otimes
u)=c_2\pi(a)u$, for some $c_1,c_2\in\C$.  Thus the formula in Proposition 4.2
(ii) reduces to
$$a_{\lambda}u=\pi(a)u
+\sum_{M\ge 1}c_M((\dd+\lambda)^{M}-\dd^M)\pi(a)u,\quad c_M\in\C.$$
By Remark 4.1 such extensions are equivalent to the trivial extension. $\qed$

For Section 5 we need to understand formula (4.6) better.  In fact we will need

\proclaim{Lemma 4.3.} Let $(\pi,U)$ and $(\rho,V)$ be non-trivial non-isomorphic finite-dimensional irreducible representations of $sl_2$ of dimensions $m+1$ and $n+1$, respectively.  Let $\varphi^{1,1}\rightarrow\varphi^{2,0}$ be the linear map defined by (4.6).  Then the map $\varphi^{1,1}\rightarrow\varphi^{2,0}$ is a bijection. Furthermore $\varphi^{2,0}=c\varphi^{1,1}$, where $c={{m+4}\over 2}$ if $n=m+2$, and $c={{2-m}\over 2}$ if $n=m-2$.

\proof Note that in order for $\varphi^{1,1}$ and $\varphi^{2,0}$ to be non-zero and $U\not\cong V$, we must have $n=m+2$ or $n=m-2$.

Since $\varphi^{1,1}$ is a $\G$-module homomorphism, we may write (4.6) as
$$\varphi^{2,0}_{[a,b]}=\varphi^{1,1}_b\pi(a)-\varphi^{1,1}_a\pi(b)+2\varphi^{1,1}_{[a,b]}.\eqno{(4.7)}$$
Thus for the basis $e,f,h$ of $sl_2$ we obtain from (4.7)
$$\eqalignno{2\varphi^{2,0}_e=&\varphi^{1,1}_e\pi(h)-\varphi^{1,1}_h\pi(e)+4\varphi^{1,1}_e,\cr
-2\varphi^{2,0}_f=&\varphi^{1,1}_f\pi(h)-\varphi^{1,1}_h\pi(f)-4\varphi^{1,1}_f,&{(4.8)}\cr
\varphi^{2,0}_h=&\varphi^{1,1}_f\pi(e)-\varphi^{1,1}_e\pi(f)+2\varphi^{1,1}_h.\cr}$$
Suppose that $\varphi^{2,0}=0$.  We will show that $\varphi^{1,1}=0$.  For this  let $u_0$ be a highest vector of the $sl_2$-module $U$ of highest weight $m$.  Then applying (4.8) to $u_0$ (using $\varphi^{2,0}=0$) we get
$$\varphi^{1,1}_e u_0=0,\quad (m-2)\varphi^{1,1}_f u_0=\rho(f)\varphi^{1,1}_h u_0,\quad\varphi^{1,1}_h u_0=0.$$
Thus $(m-2)\varphi^{1,1}_f u_0=0$.  Hence if $m\not=2$, then $\varphi^{1,1}_a u_0=0$, for all $a\in\G$, and so $\varphi^{1,1}=0$.

Now suppose that $m=2$.  If $\varphi^{1,1}_f u_0\not=0$, then $\varphi^{1,1}_f u_0$ has maximal weight in $V$ and hence must be a non-zero scalar multiple of the highest weight vector in $V$.  But $\varphi^{1,1}_f u_0$ has weight $m-2=0$ and so $V$ has highest weight $0$, which is a contradiction, since $V$ is nontrivial.  Thus $\varphi^{1,1}_fu_0=0$ and so in this case $\varphi^{1,1}=0$ as well.  This proves the first claim of the lemma.

We have $\varphi^{2,0}=c\varphi^{1,1}$.  Applying (4.8) to $u_0$ we get
$$\eqalignno{2c\varphi^{1,1}_eu_0=&(m+4)\varphi^{1,1}_eu_0,\cr
-2c\varphi^{1,1}_fu_0=&(m-4)\varphi^{1,1}_fu_0-\varphi^{1,1}_h\pi(f)u_0,\cr
c\varphi^{1,1}_h u_0=&-\varphi^{1,1}_e\pi(f)u_0+2\varphi^{1,1}_h u_0,\cr}$$
from which we derive (since $\varphi^{1,1}$ is a $\G$-homomorphism) the following equations:
$$(m+4-2c)\varphi^{1,1}_eu_0=0,\quad
(m-2+2c)\varphi^{1,1}_fu_0=\rho(f)\varphi^{1,1}_hu_0,\quad
(1-c)\varphi^{1,1}_hu_0=\rho(f)\varphi^{1,1}_eu_0.$$
If neither of the numbers $(m+4-2c),(m-2+2c),(1-c)$ is $0$, then $\varphi^{1,1}=0$ and so $\varphi^{2,0}=0$.  So we can assume that $\varphi^{1,1}\not=0$ and one of the numbers $(m+4-2c),(m-2+2c),(1-c)$ is $0$.  It is easy to see, since $m>0$, that only one of the numbers can be $0$.  Thus we consider three cases:
\item{(i)} $m+4-2c=0$ and so $\varphi^{1,1}_eu_0\not=0$ and $\varphi^{1,1}_f=\varphi^{1,1}_h=0$.
\item{(ii)} $m-2+2c=0$ and so $\varphi^{1,1}_fu_0\not=0$ and $\varphi^{1,1}_e=\varphi^{1,1}_h=0$.
\item{(iii)} $1-c=0$ and so $\varphi^{1,1}_hu_0\not =0$ and $\varphi^{1,1}_e=0$.

Now case (iii) cannot occur, because by assumption $U\not\cong V$.  case (i) means $n=m+2$, while case (ii) means $n=m-2$.  This completes the proof.  $\qed$

For simple Lie algebras of higher ranks, we can weaken the assumptions of Lemma
4.2 and still prove the same result.  We will need this in order to understand
(4.6).

\proclaim{Lemma 4.4.} Let rank $\G\ge 2$ and let $(\pi,U)$ and $({\rho},V)$ be non-trivial irreducible $\G$-modules.  Let $T:\G\otimes U\rightarrow V$ be a
$\G$-homomorphism such ${\rho}(a)T(b\otimes u)=T(a\otimes\pi(b)u)$, for all
$a,b\in\G$ with $[a,b]=0$.
\item{(i)} If $(\pi,U)\not\cong ({\rho},V)$, then $T=0$.
\item{(ii)} If $(\pi,U)\cong ({\rho},V)$, then there exists $c\in\C$ such that
$T(a\otimes u)=c\pi(a)u$, for all $a\in\G,u\in U$.

Before proving Lemma 4.4, let us see how it can be applied to
solve our extension problem in the higher rank case.

\proclaim{Proposition 4.3.} Let rank $\G\ge 2$
and let $(\pi,U)$ and
$({\rho},V)$ be non-trivial irreducible $\G$-modules.  Let $\varphi^{1,1},\varphi^{2,0}:\G\otimes
U\rightarrow V$ be two $\G$-module
homomorphisms satisfying (4.6).  Then
\item{(i)} If $(\pi,U)\not\cong({\rho},V)$, then $\varphi^{1,1}=\varphi^{2,0}=0$.
\item{(ii)} If $(\pi,U)\cong({\rho},V)$, then there exists $c\in\C$ such that
$\varphi^{1,1}_a=\varphi^{2,0}_a=c\pi(a)$, $\forall a\in\G$.

\proof Let $a,b\in\G$ such that $[a,b]=0$.  We have
${\rho}(a)\varphi^{1,1}_b={\rho}(b)\varphi^{1,1}_a$ by assumption.  Since
$\varphi^{1,1}$ is a $\G$-module homomorphism we also have
${\rho}(a)\varphi^{1,1}_b=\varphi^{1,1}_b\pi(a)$.  Hence
${\rho}(b)\varphi^{1,1}_a=\varphi^{1,1}_b\pi(a)$. So by Lemma 4.4 and (4.6) Corollary
follows. $\qed$

Combining Proposition 4.2 and Proposition 4.3 we obtain immediately

\proclaim{Theorem 4.2.} Let rank $\G\ge 2$ and let $(\pi,U)$ and $({\rho},V)$ be finite-dimensional non-trivial irreducible
$\G$-modules. Then every non-trivial extension of $R(\tilde{\G})$-modules of the form (4.1) is as follows ($a\in\G$, $u\in U$, $v\in V$):
\item{(i)} If $(\pi,U)\not\cong({\rho},V)$, then
$$a_{\lambda}u=\pi(a)u+\lambda\varphi^{1,0}(a\otimes u),\quad a_{\lambda}v=\rho(a)v,$$
where $\varphi^{1,0}:\G\otimes U\rightarrow V$ is an arbitrary non-zero homomorphism of
$\G$-modules. Two extensions, depending on $\varphi^{1,0}$ and $\tilde{\varphi}^{1,0}$, respectively, are equivalent if and only if $\varphi^{1,0}$ and $\tilde{\varphi}^{1,0}$ coincide.
\item{(ii)} If $(\pi,U)\cong({\rho},V)$, then
$$a_{\lambda}u=\pi(a)u+\lambda\varphi^{1,0}(a\otimes u),\quad a_{\lambda}v=\rho(a)v.$$
Furthermore two extensions, depending on $\varphi^{1,0}$ and
$\tilde{\varphi}^{1,0}$, respectively, are equivalent if and only if
$\varphi^{1,0}(a\otimes u)-\tilde{\varphi}^{1,0}(a\otimes u)=c\pi(a)u$, for
some $c\in\C$.

We now give a proof of Lemma 4.4.

\noindent{\it Proof of Lemma 4.4.} We will continue to use the notation above.
Let
${\alpha}_i$ be any simple root.  Since rank $\G\ge 2$ we
can choose an $h\in\H$ such that ${\alpha}_i(h)=0$ and $\Lambda(h)\not=0$. Since $[h,e_i]=0$, we have ${\rho}(h)T_{e_i}u_0=T_h\pi(e_i)u_0=0$.  Now $T_{e_i}u_0$ is a vector in $V$
of weight $\Lambda+{\alpha}_i$. Thus we have
$(\Lambda+{\alpha}_i)(h)T_{e_i}u_0=\Lambda(h)T_{e_i}u_0=0$, and so
$T_{e_i}u_0=0$.  Since $T$ is a $\G$-homomorphism, it follows that
$T(e_{{\alpha}}\otimes u_0)=0$, for ${\alpha}\in\Delta_+$.

Let $h\in\H$ and ${\alpha}_i$ be a simple root.  We have
${\rho}(e_i)T_hu_0-T_h\pi(e_i)u_0=T_{[e_i,h]}u_0=-{\alpha}_i(h)T_{e_i}u_0=0.$
Hence $T_hu_0$ is a scalar multiple of $v_0$, i.e.~there exists $c\in\H^*$ such
that
$$T_hu_0=c(h)v_0,\quad\forall h\in\H.\eqno{(4.9)}$$
So if $h,h'\in\H$ then
${\rho}(h')T_hu_0=T_h\pi(h')u_0$, which implies that
$$c(h)\Lambda'(h')v_0=c(h)\Lambda(h')v_0.\eqno{(4.10)}$$

\noindent{\csc Case 1.} Assume that $\Lambda\not=\Lambda'$. In this case (4.10)
implies that
$c(h)=0$, hence
$$T_hu_0=0,\quad\forall h\in\H\ {\rm and}\ \Lambda\not=\Lambda'.\eqno{(4.11)}$$

Let ${\alpha}\in\Delta_+$ and suppose that $\Lambda$ is not a scalar multiple of ${\alpha}$.  Then we may choose
$h\in\H$ such that
${\alpha}(h)=0$ and $\Lambda(h)\not=0$.  By (4.11)
$0={\rho}(e_{-{\alpha}})T_hu_0=T_{e_{-{\alpha}}}\pi(h)u_0=
\Lambda(h)T_{e_{-{\alpha}}}u_0$. We conclude that $T_{e_{-{\alpha}}}u_0=0$.

Now suppose that $\Lambda=c'\theta_0$, where $\theta_0\in\Delta_+$.
If $\theta_0$ is not the highest root, then there exist ${\alpha},{\beta}\in\Delta_+$ such that $\theta_0={\alpha}-{\beta}$ so that $e_{-\theta_0}=[e_{{\beta}},e_{-{\alpha}}]$.  But $T_{e_{-\theta_0}}u_0={\rho}(e_{{\beta}})T_{e_{-{\alpha}}}u_0- T_{e_{-{\alpha}}}\pi(e_{{\beta}})u_0=0$.

On the other hand if $\theta_0$ is the highest root, then $\theta_0$ cannot be a simple root.  Thus
$${\rho}(e_i)T_{e_{-\theta_0}}u_0-T_{e_{-\theta_0}}\pi(e_i)u_0=T_{[e_i,e_{-\theta_0}]}u_0=T_{e_{{\alpha}-\theta_0}}u_0=0,$$ 
hence ${\rho}(e_i)T_{e_{-\theta_0}}u_0=0$, $\forall i$.  This implies that $T_{e_{-\theta_0}}=c_0v_0$, for some $c_0\in\C$.  We need to show that $c_0=0$.  Now $[e_{-\theta_0},f_i]=0$, $\forall i$.  Therefore
$0={\rho}(e_{-\theta_0})T_{f_{i}}u_0=T_{e_{-\theta_0}}\pi(f_i)u_0= {\rho}(f_i)T_{e_{-\theta_0}}u_0$.  This gives $c_0{\rho}(f_i)v_0=0$, $\forall i$.  But $V$ is non-trivial, hence $c_0=0$.

\noindent{\csc Case 2.} Assume that $\Lambda=\Lambda'$ and let $h,h'\in\H$.  By assumption $\pi(h)T_{h'}u_0=T_h\pi(h')u_0$, hence by (4.9)
$\Lambda(h)c(h')=\Lambda(h')c(h)$, for all $h,h'\in\H$.
Since
$\Lambda\not=0$, we may choose $h_0\in\H$ such that $\Lambda(h_0)\not=0$.  
Putting $c_0={{c(h_0)}\over{\Lambda(h_0)}}$ we have $c(h)=c_0\Lambda(h)$,
$\forall h\in\H$, thus by (4.9)
$$T_hu_0=c_0\Lambda(h)u_0,\quad\forall h\in\H.\eqno{(4.12)}$$

Given ${\alpha}\in\Delta_+$ and suppose that $\Lambda$ is not a scalar multiple of ${\alpha}$.  We may choose $h\in\H$ such that ${\alpha}(h)=0$ and $\Lambda(h)\not=0$.  Then by hypothesis $\pi(e_{-{\alpha}})T_hu_0=T_{e_{-{\alpha}}}\pi(h)u_0$, which by (4.12) is equivalent to $c_0\pi(e_{-{\alpha}})\Lambda(h)u_0=\Lambda(h)T_{e_{-{\alpha}}}u_0$ and so
$$T_{e_{-{\alpha}}}u_0=c_0\pi(e_{-{\alpha}})u_0.\eqno{(4.13)}$$

Now suppose that $\Lambda=c'\theta_0$, where $\theta_0\in\Delta_+$. It suffices to prove that (4.13) when ${\alpha}=\theta_0$.
Since $T$ is a $\G$-module homomorphism, we have
$$\pi(e_i)T_{e_{-\theta_0}}u_0=T_{e_{-\theta_0}}\pi(e_i)u_0+T_{[e_i,e_{-\theta_0}]}u_0=c_0\pi([e_i,e_{-\theta_0}])u_0.$$
Replacing $\pi([e_i,e_{-\theta_0}])u_0$ by $\pi(e_i)\pi(e_{-\theta_0})u_0$, we conclude that $T_{e_{-\theta_0}}u_0-c_0\pi(e_{-\theta_0})u_0$, if nonzero, is a scalar multiple of $u_0$.  However, since it has weight $\Lambda-\theta_0$ and $u_0$ has weight $\Lambda$, we must have $T_{e_{-\theta_0}}u_0=c_0\pi(e_{-\theta_0})u_0$.  This concludes the proof of Lemma 4.4. $\qed$

We now present explicit formulas for the corresponding (non-trivial) extensions of current conformal modules. 
Recall that $\tilde{\G}=\G\otimes\C[t,t^{-1}]$ acts on the conformal module $M^c(V)=V\otimes_{\c}\C[t,t^{-1}]$, where $(\rho,V)$ is a non-trivial irreducible finite-dimensional representation of $\G$, as
$$(a\otimes f(t))(v\otimes g(t))=\rho(a)v\otimes f(t)g(t).\eqno{(4.14)}$$

\proclaim{Theorem 4.3.} All non-trivial extensions of conformal modules over $\tilde{\G}$ of the form
$$0\longrightarrow M^c(V)\longrightarrow E^c\longrightarrow M^c(U)\longrightarrow 0,$$
where $(\pi,U)$ and $(\rho,V)$ are non-trivial irreducible representations of $\G$,
can be constructed as follows. As a vector space $E^c=(U\oplus V)\otimes_{\c}\C[t,t^{-1}]$.  $\tilde{\G}$ acts on $M^c(V)$ as in (4.14), while on $M^c(U)$ its action is given as follows ($f(t),g(t)\in\C[t,t^{-1}]$, $u\in U$, $v\in V$, $a\in\G$):
$$\eqalign{(i)\ (a\otimes f(t))(u\otimes g(t))&=\pi(a)u\otimes f(t)g(t)+\varphi^{1,0}(a\otimes u)\otimes f'(t)g(t)\cr
&+(\varphi^{2,0}
-\varphi^{1,1})(a\otimes u)\otimes f''(t)g(t)
-\varphi^{1,1}(a\otimes u)\otimes f'(t)g'(t),\qquad\cr}$$
if $\G\cong sl_2$ and $V\not\cong U$.
\item{(ii)} $(a\otimes f(t))(u\otimes g(t))=\pi(a)u\otimes f(t)g(t)+\varphi^{1,0}(a\otimes u)\otimes f'(t)g(t),$
if $\G\not\cong sl_2$.\hfill\break
Here $\varphi^{1,0},\varphi^{1,1}:\G\otimes U\rightarrow V$ are $\G$-module homomorphisms, and $\varphi^{2,0}:\G\otimes U\rightarrow V$ is a $\G$-module homomorphism determined by (4.6).

\noindent{\it Remark 4.3.} In the case $\G\not\cong sl_2$, the space ${\rm Ext}(M^c(U),M^c(V))$ is canonically isomorphic to ${\rm Hom}_{\g}(\G\otimes_{\c}U,V)$ if $U\not\cong V$ and to ${\rm Hom}_{\g}(\G\otimes U,V)/\C\pi$ if $U\cong V$.  In the case $\G\cong sl_2$ and $U\not\cong V$, ${\rm Ext}(M^c(U),M^c(V))$ is canonically isomorphic to the direct sum of two copies of ${\rm Hom}_{\g}(\G\otimes_{\c}U,V)$.

The remainder of this section is devoted to the study of extensions of $R(\tilde{\G})$-modules that we will need for the extension problem of modules over $R(\V)\lsemi R(\tilde{\G})$ in the next section.  Namely we will consider extensions of the form (4.1) but now with the assumption that
$(\rho,V)$ or $(\pi,U)$ is the trivial ($1$-dimensional) $\G$-module.

First suppose that $V=\C v$ is the trivial $\G$-module. We have for $a\in\G$ and $u\in U$:
$$a_{\lambda}u=\pi(a)u+\sum_{i,k\ge 0}\lambda^i\dd^k\varphi_a^{i,k}(u),\quad a_{\lambda}v=0,$$
where $\varphi_a^{i,k}:U\rightarrow V$ are linear maps.  Applying both sides of (1.16) to $u$ we have
$$\sum_{i,k\ge 0} (\lambda+\mu)^{i}\dd^k\varphi^{i,k}_{[a,b]}(u)=\sum_{i,k\ge 0}\lambda^i\dd^k\varphi^{i,k}(\pi(b)u)
-\sum_{i,k\ge 0}\mu^i\dd^k\varphi^{i,k}(\pi(a)u).$$
Comparing the coefficients of $\lambda^i\mu^j\dd^k$ for $i+j=M$ and $i,j\ge 1$ we get for all $k$
$$\varphi_{[a,b]}^{M,k}=0,\quad\forall M\ge 2.$$
Subtracting a coboundary as in the proof of Lemma 2.1, if necessary, we may again assume that $\varphi^{0,k}=0$, for all $k$.  Now we look at the coefficients of $\lambda^i\dd^k$ and conclude that
$$\varphi^{i,k}_{[a,b]}(u)=\varphi^{i,k}_{a}(\pi(b)u),$$ 
which implies that $\varphi^{i,k}=0$, unless $U\cong \G$.  Thus we have

\proclaim{Proposition 4.4.} Non-trivial extensions of $R(\tilde{\G})$-modules of the form (4.1) with $V{=}\C v$ being the trivial and $U$ an irreducible (possibly trivial) $\G$-module exist if and only if $U\cong \G$.  Identifying $U$ with $\G$ these extension are given by ($a,u\in\G$):
$$a_{\lambda}u=[a,u]+\lambda f(\dd)(a|u),\quad a_{\lambda}v=0,$$
where $f(\dd)$ is a non-zero polynomial in $\dd$ and $(\cdot|\cdot)$ is an invariant symmetric non-degenerate bilinear form on $\G$.

Now suppose that $V$ is non-trivial and $U=\C u$ is the trivial $\G$-module. We have
$$a_{\lambda}u=\sum_{i,j\ge 0}\lambda^i\dd^j\varphi^{i,j}_a(u),\quad a_{\lambda}v=\rho(a)v,$$
where $\varphi^{i,j}_a:U\rightarrow V$ are linear maps.  Both sides of (1.16) applied to $u$ gives
$$\sum_{i,j}(\lambda+\mu)^i\dd^j\varphi^{i,j}_{[a,b]}(u)=\sum_{i,j}\mu^i(\dd+\lambda)^j\rho(a)\varphi^{i,j}_b(u)-\sum_{i,j}\lambda^i(\dd+\mu)^j\rho(b)\varphi^{i,j}_a(u).\eqno{(4.15)}$$
As before we may assume that the summation above is over $i\ge 1$ and $j\ge 0$.  Looking at the coefficients of $\mu^i\dd^j$ and we get
$\rho(a)\varphi^{i,j}_b(u)=\varphi^{i,j}_{[a,b]}(u),$
which means that $\varphi^{i,j}:\G\rightarrow V$ given by $\varphi^{i,j}(a)=\varphi^{i,j}_a(u)$ is a $\G$-module homomorphism.  Thus we may assume that $V\cong\G$.  Identifying $V$ with $\G$ we may assume that for all $a\in\G$ one has $\varphi^{i,j}_a(u)=c_{ij}a$ for some $c_{ij}\in\C$, and hence $\rho(a)\varphi^{i,j}_b(u)=c_{ij}[a,b]$.  Plugging this back in (4.15) we get
$$\sum_{i,j}(\lambda+\mu)^i\dd^j c_{ij}[a,b]=\sum_{i,j}\mu^i(\dd+\lambda)^jc_{ij}[a,b]-\sum_{i,j}\lambda^i(\dd+\mu)^jc_{ij}[a,b],$$
which is equivalent to
$$f(\lambda+\mu,\dd)=f(\mu,\dd+\lambda)+f(\lambda,\dd+\mu),\eqno{(4.16)}$$
where $f(\lambda,\dd)=\sum_{i,j}c_{ij}\lambda^i\dd^j$.  So solving the extension problem is equivalent to solving the functional equation (4.16).  We have the following

\proclaim{Lemma 4.4.} All homogeneous solutions of (4.16) are scalar multiples of $\lambda(\dd+\lambda)^k$, where $k\in\Z_+$.

\proof We use induction on the degree of a homogeneous solution of (4.7).  If the degree is $1$ or $2$ this is easily checked.  Hence assume that the degree is $n>1$.  Now differentiating (4.16) with respect to $\dd$ we obtain
$$f_{\dd}(\lambda+\mu,\dd)=f_{\dd}(\mu,\dd+\lambda)+f_{\dd}(\lambda,\dd+\mu),$$
where $f_{\dd}(\lambda,\dd)$ means differentiation with respect to $\dd$.  But this equation means that $f_{\dd}(\lambda,\dd)$ is a solution of (4.16), but of lower degree.  Thus by induction it $f_{\dd}(\lambda,\dd)$ a scalar multiple of $\lambda(\dd+\lambda)^{k-1}$, where $k<n$.  Hence
$f(\lambda,\dd)$ is a scalar multiple of $\lambda(\dd+\lambda)^{k} + g(\lambda)$, where $g(\lambda)$ is scalar multiple of $\lambda^{n}$.  Plugging this into (4.16) again, we see that $g(\lambda+\mu)=g(\lambda)+g(\mu)$.  Since $n>1$, $g(\lambda)=0$, proving the lemma. $\qed$

\proclaim{Proposition 4.5.} Non-trivial extensions of $R(\tilde{\G})$-modules of the form (4.1) with $U{=}\C u$ being the trivial and $V$ an irreducible $\G$-module exist if and only if $V\cong\G$.  Identifying $V$ with $\G$ these extensions are given by ($a,v\in\G$):
$$a_{\lambda}u=\lambda f(\dd+\lambda)a,\quad a_{\lambda}v=[a,v],\qquad
0\not=f(x)\in\C[x].$$

\beginsection{5. Extensions of conformal $\V\lsemi\tilde{\G}$-modules}

Throughout this section $\G$ denotes a finite-dimensional simple Lie algebra. Let $R(\V)\lsemi R(\tilde{\G})$ be the semidirect sum of the Virasoro and the current conformal algebra.  Any finite irreducible module over $R(\V)\lsemi R(\tilde{\G})$ is $M(\alpha,\Delta,U)$, with actions given by (1.21), or else it is the $1$-dimensional (over $\C$) module $\C c_{\beta}$ with actions $L_{(n)}c_{\beta}=a_{(n)}c_{\beta}=0$ for all $n\in\Z_+$, and $\dd c_{\beta}=\beta c_{\beta}$.  In this section we will consider extensions of $R(\V\lsemi\tilde{\G})$-modules of the form
$$0\longrightarrow M(\alpha,\bar{\Delta},V)\longrightarrow E\longrightarrow M(\beta,\Delta,U)\longrightarrow 0,\eqno{(5.1)}$$
where $(\rho,V)$ and $(\pi,U)$ are finite-dimensional irreducible $\G$-modules. 

\noindent {\it Remark 5.1.} We do not need to consider the case when $U$ and $V$ are both the trivial $\G$-module, because due to Proposition 4.4 (or 4.5) this case reduces to the extension problem we have studied in Section 3.

We first consider the case when $U$ and $V$ are both non-trivial $\G$-modules.  In this case we can use the results in Section 4 and assume that ($a\in\G$, $u\in U$ and $v\in V$):
$$\eqalignno{a_{\lambda}u&=\pi(a)u+\lambda\varphi^{1,0}(a\otimes u)+2\lambda\dd\varphi^{1,1}(a\otimes u) +\lambda^2\varphi^{2,0}(a\otimes u),\cr
a_{\lambda}v&=\rho(a)v,\quad
L_{\lambda}u=(\dd+\beta+\Delta\lambda)u+\sum_{j,k\ge 0}\lambda^j\dd^k\psi_{jk}u,\quad
L_{\lambda}v=(\dd+\alpha+\bar{\Delta}\lambda)v,\cr}$$
where $\varphi^{i,j}:\G\otimes U\rightarrow V$ are $\G$-module homomorphisms and $\psi_{jk}:U\rightarrow V$ are linear maps.  Furthermore $\varphi^{2,0}([a,b]\otimes u)=\rho(a)\varphi^{1,1}(b\otimes u)-\rho(b)\varphi^{1,1}(a\otimes u)$ and $\varphi^{1,1}$ and $\varphi^{2,0}$ are non-zero only when $\G\cong sl_2$ and $U\not\cong V$.

Applying both sides of (1.18) to $u$ we obtain:
$$\eqalignno{&-\lambda\mu\varphi^{1,0}_a u-2\mu(\lambda+\mu)\dd\varphi^{1,1}_a u-\mu(\lambda+\mu)^2\varphi^{2,0}_a u=\mu(\alpha-\beta+(\bar{\Delta}-\Delta)\lambda)\varphi^{1,0}_a u\cr
&+2\mu((\dd+\alpha+\lambda)(\dd+\bar{\Delta}\lambda)-\dd(\dd+\beta+\mu+\Delta\lambda))\varphi^{1,1}_a u+\mu^2(\alpha-\beta-\mu-(\bar{\Delta}-\Delta)\lambda)\varphi^{2,0}_a u\cr
&+\sum_{j,k\ge 0}\lambda^j\dd^k\psi_{jk}\pi(a)u-\sum_{j,k\ge 0}\lambda^j(\dd+\mu)^k\rho(a)\psi_{jk}u.&(5.2)\cr}$$
Setting $\mu=0$ in (5.2) and comparing the coefficients of $\lambda^j\dd^k$  we obtain for all $j,k\ge 0$
$$\psi_{jk}\pi(a)u=\rho(a)\psi_{jk}u,\quad\forall a\in\G,u\in U.$$
Thus $\psi_{jk}$ is a $\G$-module homomorphism.

Assume that $U\not\cong V$.  If $\G\not\cong sl_2$, then $\varphi^{1,1}=\varphi^{2,0}=0$ and so (5.2), combined with the fact that $\psi_{jk}=0$, reduces to
$$((\Delta-\bar{\Delta}-1)\lambda +\beta-\alpha)\varphi^{1,0}_a u=0.$$
If $\varphi^{1,0}=0$, then $E$ decomposes.  Thus we may assume that $\varphi^{1,0}\not=0$.  In this case we get $\alpha=\beta$ and $\Delta-\bar{\Delta}-1=0$.  This gives

\proclaim{Proposition 5.1.} Non-trivial extensions of $R(\V\lsemi\tilde{\G})$-modules of the form (5.1) in the case $\G\not\cong sl_2$, $(\pi,U)\not\cong(\rho,V)$ and $U,V$ non-trivial irreducible $\G$-modules exist only if $\alpha=\beta$ and $\Delta=\bar{\Delta}+1$.  The extensions are given by:
$$\eqalign{a_{\lambda}u=&\pi(a)u+\lambda\varphi^{1,0}(a\otimes u),\qquad
a_{\lambda}v=\rho(a)v,\cr
L_{\lambda}u=&(\dd+\alpha+\Delta\lambda)u,\hskip 0.63in
L_{\lambda}v=(\dd+\alpha+\bar{\Delta}\lambda)v,\cr
}$$
where $u\in U$, $v\in V$, $a\in\G$ and $\varphi^{1,0}:\G\otimes U\rightarrow V$ is a non-trivial $\G$-module homomorphism.

Now suppose that $\G\cong sl_2$ and $U\not\cong V$.  In this case (5.2) together with $\psi_{jk}=0$ gives:
$$\eqalignno{&-\lambda\varphi^{1,0}_a u-2(\lambda+\mu)\dd\varphi^{1,1}_a u-(\lambda+\mu)^2\varphi^{2,0}_a u=(\alpha-\beta+(\bar{\Delta}-\Delta)\lambda)\varphi^{1,0}_a u&{(5.3)}\cr
&+2((\dd+\lambda)(\dd+\alpha+\bar{\Delta}\lambda)-\dd(\dd+\beta+\mu+\Delta\lambda))\varphi^{1,1}_a u+\mu(\alpha-\beta-\mu-(\bar{\Delta}-\Delta)\lambda)\varphi^{2,0}_a u\cr}$$
Putting $\lambda=0$, (5.3) reduces to
$(\alpha-\beta)(\varphi^{1,0}_a u+2\dd\varphi^{1,1}_a u+\mu\varphi^{2,0}_a u)=0$. Hence if $\alpha-\beta\not=0$, then $\varphi^{i,j}=0$ for all $i,j$ and so $E$ decomposes.  Thus we may assume that $\alpha-\beta=0$.  In this case (5.3) reduces to
$$\lambda(\Delta-\bar{\Delta}-1)\varphi^{1,0}_a u-\lambda^2(2\bar{\Delta}\varphi^{1,1}_a u+\varphi^{2,0}_a u)+2\lambda\dd(\Delta-\bar{\Delta}-2-\alpha)\varphi^{1,1}_a u+\lambda\mu(\Delta-\bar{\Delta}-2)\varphi^{2,0}_a u=0.$$
From this and Lemma 4.3 we see that $\varphi^{i,j}=0$ and hence the representation decomposes, if $\Delta-\bar{\Delta}\not=1,2$.  In the case when $\Delta-\bar{\Delta}=1$, $\varphi^{2,0}=0$ and hence $\varphi^{1,1}=0$ by Lemma 4.3.  On the other hand when $\Delta-\bar{\Delta}=2$ and $\alpha=0$ we have $\varphi^{1,0}=0$ and  $\varphi^{2,0}=-2\bar{\Delta}\varphi^{1,1}$.  Now if $\Delta-\bar{\Delta}=2$  and $\alpha\not=0$, then $\varphi^{1,0}=0$.  But also $\varphi^{1,1}=0$, which gives $\varphi^{2,0}=0$ by Lemma 4.3.  Hence in this case the module decomposes.  This discussion, together with Lemma 4.3, gives

\proclaim{Proposition 5.2.} Let $\G\cong sl_2$. Non-trivial extensions of $R(\V\lsemi\tilde{\G})$-modules of the form (5.1) in the case $(\pi,U)\not\cong(\rho,V)$ with $U,V$ non-trivial irreducible $\G$-modules are as follows:
\item{(i)} $\alpha=\beta$ and $\Delta-\bar{\Delta}=1$, with extensions given as in Proposition 5.1.
\item{(ii)} $\alpha=\beta=0$ and $\Delta-\bar{\Delta}=2$, with extensions given by ($u\in U$, $v\in V$, $a\in\G$):
$$\eqalign{a_{\lambda}u=&\pi(a)u+2\lambda\dd\varphi^{1,1}(a\otimes u)-2(\Delta-2)\lambda^2\varphi^{1,1}(a\otimes u),\cr
a_{\lambda}v=&\rho(a)v,\quad
L_{\lambda}u=(\dd+\Delta\lambda)u,\quad
L_{\lambda}v=(\dd+(\Delta-2)\lambda)v,\cr
}$$
where $\varphi^{1,1}:\G\otimes U\rightarrow V$ is a non-trivial $\G$-module homomorphism and $\Delta={{5-{\rm dim}U}\over 4}$, if ${\rm dim}V={\rm dim}U+2$, and $\Delta={{5+{\rm dim}U}\over 4}$, if ${\rm dim}V={\rm dim}U-2$.

This settles the case $U\not\cong V$ and $U$ and $V$ are non-trivial $\G$-modules.  Next we consider the case $U\cong V$ and $U$ and $V$ are non-trivial $\G$-modules.  In this case we have by results of Section 4 that $\varphi^{1,1}=\varphi^{2,0}=0$.  Also in this case $\psi_{jk}(u)=c_{jk}u$, where $c_{jk}\in\C$.  Thus identifying $\rho$ with $\pi$, (5.2) reduces to
$$-\lambda\mu\varphi^{1,0}_a u=\mu(\alpha-\beta+(\bar{\Delta}-\Delta)\lambda)\varphi^{1,0}_a u+\sum_{j,k\ge 0}c_{jk}\lambda^j(\dd^k-(\dd+\mu)^k)\pi(a)u.\eqno{(5.4)}$$
We let $\dd=0$ and get
$$\mu(\alpha-\beta+(\bar{\Delta}-\Delta+1)\lambda)\varphi^{1,0}_a u-\sum_{j\ge 0,k\ge 1}c_{jk}\lambda^j\mu^k\pi(a)u=0.\eqno{(5.5)}$$
If $\varphi^{1,0}$ is a scalar multiple of $\pi$, then it is decomposable and we may assume that $\varphi^{1,0}=0$. 
If  
$\varphi^{1,0}=0$, then $c_{jk}=0$ for $k\ge 1$.  Thus  
$L_{\lambda}u=(\partial +\beta+\Delta\lambda)u+f(\lambda)u'$, where  
$u\rightarrow u'$ is a $\G$-module isomorphism from $U$ to $V$.  Since  
$f(\lambda)$ must satisfy $(3.2)$, we get from the list of Theorem 3.1 the  
following possibilities:
\item{(I1)} $f(\lambda)=a_0+a_1\lambda$ in the case $\Delta=\bar{\Delta}$,  
$(a_0,a_1)\not=(0,0)$,
\item{(I2)} $f(\lambda)=a\lambda^2$ in the case $\Delta=1$ and  
$\bar{\Delta}=0$, $a\not=0$,
\item{(I3)} $f(\lambda)=a\lambda^3$ in the case $\Delta-\bar{\Delta}=2$,  
$a\not=0$.

\noindent
Thus we may assume that $\varphi^{1,0}$ is not a scalar multiple of $\pi$. In this case (5.5) implies that $\alpha-\beta=0$, $\Delta=\bar{\Delta}+1$ and $c_{jk}=0$ for $j\ge 0$ and $k\ge 1$.  Thus setting $\sum_{jk}c_{jk}\lambda^j\dd^k=f(\lambda)$, we have (where $u\rightarrow\bar{u}$ is a $\G$-module isomorphism between $U$ and $V$)
$$L_{\lambda}u=(\dd+\alpha+\Delta\lambda)u+f(\lambda)\bar{u},\quad L_{\lambda}v=(\dd+\alpha+(\Delta-1)\lambda)v.$$
Applying both sides of (1.17) to $u$ we obtain
$$(\lambda-\mu)f(\lambda+\mu)=(\lambda+\mu)(f(\lambda)-f(\mu)).$$
The only solution to this equation is $f(\lambda)=c\lambda$, where $c\in\C$.  However this polynomial corresponds to the trivial extension.  Namely let $E=\C[\dd]U\oplus\C[\dd]\bar{U}$ and $u'=u+\bar{u}$.  Then $L_{\lambda}u'$ gives rise to this polynomial.  Thus we have

\proclaim{Proposition 5.3.} Non-trivial extensions of $R(\V\lsemi\tilde{\G})$-modules of the form (5.1) in the case $(\pi,U)\cong(\rho,V)$ and $U,V$ non-trivial irreducible $\G$-modules are either of types (I1), (I2) and (I3) above or else exist only if $\G\not\cong sl_2$, $\alpha=\beta$ and $\Delta=\bar{\Delta}+1$.  They are given by
$$\eqalign{a_{\lambda}u=&\pi(a)u+\lambda\varphi^{1,0}(a\otimes u),\qquad
a_{\lambda}v=\rho(a)v,\cr
L_{\lambda}u=&(\dd+\alpha+\Delta\lambda)u,\hskip 0.65in
L_{\lambda}v=(\dd+\alpha+(\Delta-1)\lambda)v,\cr
}$$
where $u\in U$, $v\in V$, $a\in\G$ and $\varphi^{1,0}:\G\otimes U\rightarrow V$ is a $\G$-module homomorphism that is not a scalar multiple of $\pi$.

In terms of conformal modules Proposition 5.1, 5.2 and 5.3 gives

\proclaim{Theorem 5.1.} All extensions of conformal modules over $\V\lsemi\tilde{\G}$ of the form
$$0\longrightarrow M^c(\alpha,\bar{\Delta},V)\longrightarrow E^c\longrightarrow M^c(\beta,\Delta,U)\longrightarrow 0,$$
where $(\pi,U)$ and $(\rho,V)$ are non-trivial irreducible $\G$-modules and $\alpha,\beta,\Delta,\bar{\Delta}\in\C$, can be constructed as follows:
Let $E^c=(U\oplus V)\otimes\C[t,t^{-1}]e^{-\alpha t}$ with completely reducible action of $\V$.  Then one has the following four cases:
\item{(i)} Let $\G\cong sl_2$ and $U\not\cong V$.  Then $\tilde{\G}$ acts on $E^c$ as in Theorem 4.3 (i) with  $\varphi^{1,0}\not=0$ and $\varphi^{1,1}=\varphi^{2,0}=0$, and $\Delta=\bar{\Delta}+1$.
\item{(ii)} Let $\G\cong sl_2$ and $U\not\cong V$.  Then $\tilde{\G}$ acts on $E^c$ as in Theorem 4.3 (i) with  $\varphi^{1,0}=0$ and $\varphi^{1,1}=\varphi^{2,0}\not=0$, and $\Delta=\bar{\Delta}+2$.  Furthermore $\alpha=0$.
\item{(iii)} Let $\G\not\cong sl_2$.  Then $\tilde{\G}$ acts on $E^c$ as in Theorem 4.3 (ii) and $\Delta=\bar{\Delta}+1$.
\item{(iv)} $\G$ arbitrary and $U\cong V$.  Then $\tilde{\G}$ acts completely reducibly on $E^c$, while $\V$ acts according to (I1), (I2) and (I3).

We next consider the case when either $U$ or $V$ in (5.1) is the trivial $\G$-module, but not both.  First suppose that $V=\C v$ is the trivial $\G$-module.  In order for $M(\alpha,\bar{\Delta},V)$ to be irreducible, $\bar{\Delta}\not=0$.
By Proposition 4.4 we may put ($a\in\G$, $u\in U$ and $v\in V$)
$$\eqalign{a_{\lambda}u=&\pi(a)u+\lambda f(\dd)(a|u)v,\hskip 1.1in
a_{\lambda}v=0,\cr
L_{\lambda}u=&(\dd+\beta+\Delta\lambda)u+\sum_{j,k\ge 0}\lambda^j\dd^k\psi_{jk}(u),\qquad
L_{\lambda}v=(\dd+\alpha+\bar{\Delta}\lambda)v,\cr
}$$
where $(\cdot|\cdot)=0$, if $U\not\cong\G$, and $(\cdot|\cdot)$ is an invariant symmetric non-degenerate bilinear form on $\G$, if $U\cong\G$, $f(\dd)$ is a polynomial in $\dd$ and $\psi_{jk}:U\rightarrow V$ are linear maps.  Applying both sides of (1.18) to $u$ we get
$$\eqalignno{-\mu(\lambda+\mu)f(\dd)(a|u)&=\mu (\dd+\alpha+\bar{\Delta}\lambda)f(\dd+\lambda)(a|u)\cr
&-\mu(\dd+\mu+\beta+\Delta\lambda)f(\dd)(a|u)+\sum_{j,k\ge 0}\lambda^j\dd^k\psi_{jk}\pi(a)u.&{(5.6)}\cr}$$
Set $\mu=0$ in (5.6) above and comparing the coefficients of $\lambda^j\dd^k$ we conclude that (since $(\pi,U)$ is non-trivial) $\psi_{jk}=0$, for all $j,k$.
Now if $U\not\cong \G$, then $(\cdot|\cdot)=0$, and so the representation decomposes.  Thus we may assume that $U\cong\G$.  In this case (5.6) reduces to
$$(\dd+\beta+(\Delta-1)\lambda)f(\dd)=(\dd+\alpha+\bar{\Delta}\lambda)f(\dd+\lambda).
\eqno{(5.7)}$$
Putting $\lambda=0$ in (5.7) gives $\alpha=\beta$, if $f(\dd)\not=0$.  Thus we may assume that $\alpha=\beta$. Changing the variable $\dd$ to $\dd+\alpha$ in (5.7), we may assume that in (5.7) above $\alpha=\beta=0$.  In this case we set $\dd=0$ in (5.7) and (using $\bar{\Delta}\not=0$) conclude that $f(\dd)=c$, $c\in\C$, and $\Delta=\bar{\Delta}+1$.  This leads to

\proclaim{Proposition 5.4.} Non-trivial extensions of $R(\V\lsemi\tilde{\G})$-modules of the form (5.1) in the case when $V{=}\C v$ is the trivial and $U$ is a non-trivial irreducible $\G$-module exist if and only if $\alpha=\beta$, $\Delta=\bar{\Delta}+1$ and $U\cong\G$.  Identifying $U$ with $\G$ these extensions are given by 
$$\eqalign{a_{\lambda}u=&[a,u]+c \lambda(a|u)v,\hskip 0.33in
a_{\lambda}v=0,\cr
L_{\lambda}u=&(\dd+\alpha+\Delta\lambda)u,\hskip 0.4in
L_{\lambda}v=(\dd+\alpha+\bar{\Delta}\lambda)v,\cr
}$$
where $a,u\in\G$, $c$ is a non-zero complex number and $(\cdot|\cdot)$ is a non-degenerate symmetric invariant bilinear form on $\G$.

Finally consider the case when $(\rho,V)$ is non-trivial and $U=\C u$ is the trivial $\G$-module.  In this case $\Delta\not=0$. Using Proposition 4.5 we may assume that ($a\in\G$, $v\in V$)
$$\eqalign{a_{\lambda}u=&\lambda f(\dd+\lambda)a,\hskip 1.45in
a_{\lambda}v=\rho(a)v,\cr
L_{\lambda}u=&(\dd+\beta+\Delta\lambda)u+\sum_{j,k\ge 0}\lambda^j\dd^k v_{jk},\qquad
L_{\lambda}v=(\dd+\alpha+\bar{\Delta}\lambda)v,\cr
}$$
where $v_{jk}\in V$ and $f\not=0$ only if $V\cong \G$.  Applying both sides of (1.18) to $u$ we obtain
$$\eqalignno{-\mu(\lambda+\mu)f(\dd+\lambda+\mu)a=&\mu(\dd+\alpha+\bar{\Delta}\lambda)f(\dd+\lambda+\mu)a&{(5.8)}\cr-&\mu(\dd+\beta+\mu
+\Delta\lambda)f(\dd+\mu)a-\sum_{j,k\ge 0}\lambda^j(\dd+\mu)^k\rho(a)v_{jk}.\cr}$$
Setting $\mu=0$ in (5.8) (since $V$ is non-trivial) gives $v_{jk}=0$ for all $j,k$.  Hence, if $V\not\cong\G$, then $f=0$ and the representation decomposes.  Thus we may assume that $V\cong\G$.  Then (5.8) gives
$$-(\lambda+\mu)f(\dd+\lambda+\mu)=(\dd+\alpha+\bar{\Delta}\lambda)f(\dd+\lambda+\mu)-(\dd+\beta+\mu+\Delta\lambda)f(\dd+\mu).\eqno{(5.9)}$$
Put $\lambda=0$ in (5.9) we get $(\alpha-\beta) f(\dd+\mu)=0$.  Since the representation decomposes when $f=0$, we may assume from now on that $\alpha=\beta$.  In this case we change $\dd$ to $\dd+\alpha$ in (5.9) and thus we want to solve
$$(\dd+\Delta\lambda+\mu)f(\dd+\mu)=(\dd+(\bar{\Delta}+1)\lambda+\mu)f(\dd+\mu+\lambda).\eqno{(5.10)}$$
Put $\dd=\mu=0$ in (5.10) we get $\Delta\lambda f(0)=(\bar{\Delta}+1)\lambda f(\lambda)$.  Thus if $f(0)\not=0$, then $f(\lambda)$ is a constant.  Now if $f(0)=0$ and $f(\lambda)\not=0$, then $\bar{\Delta}+1=0$.  Now put $\mu=0$ in (5.10) and we get
$(\dd+\Delta\lambda)f(\dd)=\dd f(\dd+\lambda)$.
Hence we may write $f(\dd)=\dd g(\dd)$, for some polynomial $g(\dd)$ and consequently get
$$(\dd+\Delta\lambda)\dd g(\dd)=\dd(\dd+\lambda)g(\dd+\lambda).$$
Put $\dd=-\lambda$ and we have $(-\lambda+\Delta\lambda)(-\lambda)g(-\lambda)=0$.  Thus if $f\not=0$, then $\Delta=1$.  Now we put $\Delta=1$, $\bar{\Delta}+1=0$ and $\mu=0$ in (4.6) to get
$(\dd+\lambda)f(\dd)=\dd f(\dd+\lambda)$.
It follows that $f(\dd)=\dd g(\dd)$ and hence $(\dd+\lambda)\dd g(\dd)=(\dd+\lambda)\dd g(\dd+\lambda)$.  Thus $g$ is a constant.  Therefore $f(\lambda)$ is a scalar multiple of $\lambda$.

\proclaim{Proposition 5.5.} Non-trivial extensions of $R(\V\lsemi\tilde{\G})$-modules of the form (5.1) in the case when $U{=}\C u$ is the trivial and $V$ is a non-trivial irreducible $\G$-module exist if only if $\alpha=\beta$, $V\cong\G$ and $\Delta=\bar{\Delta}+1$ or $(\Delta,\bar{\Delta})=(1,-1)$.  Identifying $V$ with $\G$ these extensions are given by:
$$\eqalign{a_{\lambda}u=&c\lambda a,\hskip 1in
a_{\lambda}v=[a,v],\cr
L_{\lambda}u=&(\dd+\alpha+\Delta\lambda)u,\qquad
L_{\lambda}v=(\dd+\alpha+\bar{\Delta}\lambda)v,\cr
}$$
where $\Delta=\bar{\Delta}+1$, $a,v\in\G$ and $c$ is a non-zero complex number, or
$$\eqalign{a_{\lambda}u=&c\lambda(\lambda+\dd+\alpha) a,\qquad
a_{\lambda}v=[a,v],\cr
L_{\lambda}u=&(\dd+\alpha+\lambda)u,\qquad
L_{\lambda}v=(\dd+\alpha-\lambda)v,\cr
}$$
where $a,v\in\G$ and $c$ is a non-zero complex number.

Translating Proposition 5.4 and 5.5 into the language of conformal modules we obtain

\proclaim{Theorem 5.2.} All non-trivial extensions of conformal modules over $\V\lsemi\tilde{\G}$ of the form
$$0\longrightarrow M^c(\alpha,\bar{\Delta},V)\longrightarrow E^c\longrightarrow M^c(\beta,\Delta,U)\longrightarrow 0,$$
where $U$ and $V$ are irreducible $\G$-modules with one of them being the trivial $\G$-module (but not both) can be constructed as follows (for $f(t)\in\C[t,t^{-1}]$, $g(t)\in\C[t,t^{-1}]e^{-\alpha t}$, $\Delta,\alpha\in\C$, $c\not=0$ and $a,b\in\G$ and $(\cdot|\cdot)$ is an invariant symmetric non-degenerate bilinear form):
\item{(i)} Suppose that $V=\C v$ is trivial.  Then $U\cong \G$ and 
$E^c=(\G\oplus \C v)\otimes\C[t,t^{-1}]e^{-\alpha t}$ with completely reducible action of $\V$. Furthermore $\Delta=\bar{\Delta}+1$, and $\tilde{\G}$ acts as:
$$\eqalignno{&(a\otimes f(t))(b\otimes g(t))=[a,b]\otimes f(t)g(t) + c(a|b)v\otimes f'(t)g(t),\cr
&(a\otimes f(t))(v\otimes g(t))=0.\cr}$$
\item{(ii)} Suppose $U=\C u$ is trivial. Then $V\cong\G$ and $E^c=(\G\oplus \C u)\otimes\C[t,t^{-1}]e^{-\alpha t}$ with completely reducible action of $\V$, and there are two possibilities:
\item{a.} $\Delta=\bar{\Delta}+1$ and $\tilde{\G}$ acts as:
$$ 
(a\otimes f(t))(u\otimes g(t))=ca\otimes f'(t)g(t),\quad
(a\otimes f(t))(b\otimes g(t))=[a,b]\otimes f(t)g(t).$$
\item{b.} $\Delta=1$, $\bar{\Delta}=-1$ and $\tilde{\G}$ acts as:
$$(a\otimes f(t))(u\otimes g(t))=ca\otimes f'(t)g'(t),\quad
(a\otimes f(t))(b\otimes g(t))=[a,b]\otimes f(t)g(t).$$

Let $\A=\C a$ be the $1$-dimensional abelian Lie algebra so that $\tilde{\A}$ is its associated current algebra.  Let $\V\lsemi\tilde{\A}$ be its semidirect sum with the Virasoro algebra.  Let $R(\V\lsemi\tilde{\A})$ denote the corresponding conformal algebra.  Using the same argument as in [1] one can prove that every finite irreducible module over $R(\V\lsemi\tilde{\A})$ is either $M(\alpha,\Delta,k){=}\C[\dd]\otimes_{\c}\C u$, with actions defined by
$$L_{\lambda}u=(\dd+\alpha+\Delta\lambda)u,\quad a_{\lambda}u=ku,$$
where $\alpha,\Delta,k\in\C$ with $(\Delta,k)\not=(0,0)$, or else it is the $1$-dimensional module $\C c_{\beta}$ with $a_{(n)}c_{\beta}=L_{(n)}c_{\beta}=0$ and $\dd c_{\beta}=\beta c_{\beta}$, $\beta\in\C$.
Of course the corresponding conformal module $M^c(\alpha,\Delta,k)$ is the $\V$-module $\C[t,t^{-1}]e^{-\alpha t}dt^{1-\Delta}$, on which $\tilde{\A}$ acts as:
$$(a\otimes f(t))g(t)dt^{1-\Delta}=kf(t)g(t)dt^{1-\Delta},$$
where $f(t)\in\C[t,t^{-1}]$ and $g(t)\in\C[t,t^{-1}]e^{-\alpha t}$.
We now consider the extension problem between finite irreducible modules over $R(\V\lsemi\tilde{\A})$.  In the following discussion we will omit the proofs, that made use of our results on extensions between irreducible $R(\V)$-modules of Section 2 and 3.

We first consider extensions of $R(\V\lsemi\tilde{\A})$-modules of the form
$$0\longrightarrow\C c_{\beta}\longrightarrow E\longrightarrow M(\alpha,\Delta,k)\longrightarrow 0.\eqno{(5.11)}$$
Identifying $E$ with $\C c_{\beta}\oplus M(\alpha,\Delta,k)$ and $M(\alpha,\Delta,k)$ with $\C[\dd] u$ we have

\proclaim{Proposition 5.6.} The following is a complete list of non-trivial extensions of $R(\V\lsemi\tilde{\A})$-modules of the form (5.11):
\item{(i)} $L_{\lambda}u=(\dd+\alpha+\lambda)u+a_2\lambda^2c_{-\alpha}$ and $a_{\lambda}u=b_1\lambda c_{-\alpha}$, $(a_2,b_1)\not=(0,0)$.
\item{(ii)} $L_{\lambda}u=(\dd+\alpha+2\lambda)u+a_3\lambda^3c_{-\alpha}$ and $a_{\lambda}u=0$, $a_3\not=0$.
\item{(iii)} $L_{\lambda}u=(\dd+\alpha)u$ and $a_{\lambda}u=ku+b_0c_{-\alpha}$, $b_0\not=0$, $k\not=0$.

\noindent{\it Remark 5.2.} The corresponding $\V\lsemi\tilde{\A}$-modules of Proposition 5.6 can be constructed as follows ($f(t)\in\C[t,t^{-1}],g(t)\in\C[t,t^{-1}]e^{-\alpha t} $): (i) is the $\V$-module $E^c=\C[t,t^{-1}]e^{-\alpha t}+\C c_{-\alpha}$ of Remark 2.1 with $\tilde{\A}$ acting as: 
$$(a\otimes f(t))g(t)=b_1({\rm Res}_t f'(t)g(t))c_{-\alpha},$$
(ii) is $E^c=\C[t,t^{-1}]e^{-\alpha t}dt^{-1}+\C c_{-\alpha}$ of Remark 2.1 with trivial $\tilde{\A}$-action, while
(iii) is the module $E^c=\C[t,t^{-1}]e^{-\alpha t}dt+\C c_{-\alpha}$ with actions
$$\eqalignno{(a\otimes f(t))g(t)dt=&kf(t)g(t)dt+b_0({\rm Res}_tf(t)g(t))c_{-\alpha},\cr
(f(t)\dt)g(t)dt=&(f(t)g(t))'dt.\cr}$$

Next we consider extensions of $R(\V\lsemi\tilde{\A})$-modules of the form
$$0\longrightarrow M(\alpha,\Delta,k)\longrightarrow E\longrightarrow \C c_{\beta}\longrightarrow 0.\eqno{(5.12)}$$

\proclaim{Proposition 5.7.} The only non-trivial extensions of $R(\V\lsemi\tilde{\A})$-modules of the form (5.12) are non-trivial extensions of the form (2.4) with trivial action of $R(\tilde{\A})$. The corresponding conformal module is the $\V$-module of Remark 2.2 with trivial $\tilde{\A}$-action.

\proclaim{Theorem 5.3.} All non-trivial extensions over $\V\lsemi\tilde{\A}$ between a $1$-dimensional (over $\C$) module and $M^c(\alpha,\Delta,k)$ are given by Remark 5.2 and Proposition 5.7.

Now we come to extensions of $R(\V\lsemi\tilde{\A})$-modules of the form
$$0\longrightarrow M(\alpha,\bar{\Delta},\bar{k})\longrightarrow E\longrightarrow M(\beta,\Delta,k)\longrightarrow 0.\eqno{(5.13)}$$
We let $M(\alpha,\bar{\Delta},\bar{k})=\C[\dd]v$ and $M(\beta,\Delta,k)=\C[\dd]u$ so that as a $\C[\dd]$-module one has $E=M(\alpha,\bar{\Delta},\bar{k})\oplus M(\beta,\Delta,k)$ and $$\eqalignno{L_{\lambda}u=&(\dd+\beta+\Delta\lambda)u+f(\dd,\lambda)v,\quad L_{\lambda}v=(\dd+\alpha+\bar{\Delta}\lambda)v,\cr
a_{\lambda}u=&ku+g(\dd,\lambda)v,\hskip 0.9in a_{\lambda}v=\bar{k}v,\cr}$$
where $f(\dd,\lambda)$ and $g(\dd,\lambda)$ are polynomials in $\dd$ and $\lambda$ determining the extension.

\proclaim{Theorem 5.4.} The following is a complete list of non-trivial extensions of $R(\V\lsemi\tilde{\A})$-modules of the form (5.13). As before we exclude the cases $(\Delta,k)=(0,0)$ and $(\bar{\Delta},\bar{k})=(0,0)$.
\item{(i)} A non-trivial extension of $R(\V)$-modules of Theorem 3.2 with trivial $R(\tilde{\A})$-action.
\item{(ii)} $L_{\lambda}u=(\dd+\alpha+\Delta\lambda)u+(a_0+a_1\lambda)v$, $L_{\lambda}v=(\dd+\alpha+\Delta\lambda)v$, $a_{\lambda}u=ku+b_0v$ and $a_{\lambda}v=kv$, $(a_0,a_1,b_0)\not=(0,0,0)$, $\alpha,k,\Delta\in\C$.
\item{(iii)} $L_{\lambda}u=(\dd+\alpha+\Delta\lambda)u$, $L_{\lambda}v=(\dd+\alpha+(\Delta-1)\lambda)v$, $a_{\lambda}u=ku+b_1\lambda v$ and $a_{\lambda}v=kv$, $b_1\not=0$, $\alpha,k\in\C$ and $\Delta\not=1$.
\item{(iii')} $L_{\lambda}u=(\dd+\alpha+\lambda)u+a_2\lambda^2 v$, $L_{\lambda}v=(\dd+\alpha)v$, $a_{\lambda}u=ku+a_2\lambda v$ and $a_{\lambda}v=kv$, $a_2\not=0$, $\alpha\in\C$, $k\not=0$.
\item{(iv)} $L_{\lambda}u=(\dd+\alpha+\Delta\lambda)u$, $L_{\lambda}v=(\dd+\alpha+(\Delta-2)\lambda)v$, $a_{\lambda}u=ku+b_2\lambda(\dd+\alpha-(\Delta-2)\lambda) v$ and $a_{\lambda}v=kv$, $b_2\not=0$, $\alpha,k,\Delta\in\C$.
\item{(iv')} $L_{\lambda}u=(\dd+\alpha+\Delta\lambda)u+a_2\lambda^2(2\dd+2\alpha+\lambda)$, $L_{\lambda}v=(\dd+\alpha+(\Delta-2)\lambda)v$, $a_{\lambda}u=a_2\lambda(\dd+\alpha-(\Delta-2)\lambda) v$ and $a_{\lambda}v=0$, $a_2\not=0$, $\alpha,\Delta\in\C$.
\item{(v)} $L_{\lambda}u=(\dd+\alpha+\lambda)u+a_3\lambda^2(\dd+\alpha)(\dd+\alpha+\lambda)$, $L_{\lambda}v=(\dd+\alpha-2\lambda)v$, $a_{\lambda}v=0$ and \hfill\break
$a_{\lambda}u=(b_3((\dd+\alpha)^2(\dd+\alpha+\lambda)-(\dd+\alpha+\lambda)^2(\dd+\alpha-2\lambda))+ a_3(\lambda(\dd+\alpha)^2+3\lambda^2(\dd+\alpha)+2\lambda^3))v,$
with $(a_3,b_3)\not=(0,0)$ and $\alpha,\Delta\in\C$.

Theorem 5.4 above describes all extension of $\V\lsemi\tilde{\A}$-modules of the form
$$0\longrightarrow M^c(\alpha,\bar{\Delta},\bar{k})\longrightarrow E^c\longrightarrow M^c(\beta,\Delta,k)\longrightarrow 0.$$
Explicit formulas for the actions of $\V$ and $\tilde{\A}$ on the corresponding conformal modules can be derived just as we have derived (3.11).  To do so we let $E^c=\C[t,t^{-1}]e^{-\alpha t}dt^{1-\bar{\Delta}}\oplus\C[t,t^{-1}]e^{-\alpha t}dt^{1-\Delta}$.  Suppose that the extension of the corresponding module over $R(\V\lsemi\tilde{\A})$ for $\alpha=0$ is given by the pair of polynomials $f(\dd,\lambda)$ and $g(\dd,\lambda)$ in Theorem 5.4.  The action of $\V$ is of course given by (3.11), while the action of $\tilde{\A}$ is
$$(a\otimes f(t))g(t)dt^{1-\Delta}=kf(t)g(t)dt^{1-\Delta}-\sum_{i,k} b_{ik}(-1)^i (\sum_{j=0}^i{i\choose j}f^{(k+i-j)}(t)g^{(j)}(t))dt^{1-\bar{\Delta}},$$
where $g(\dd,\lambda)=\sum_{i,j}b_{ij}\dd^i\lambda^j$, $f(t)\in\C[t,t^{-1}]$ and $g(t)\in\C[t,t^{-1}]e^{-\alpha t}$.

\bigskip
\bigskip
\centerline  {\bf REFERENCES }
\bigskip
\frenchspacing
\medskip

\item{1.}Cheng, S.-J.; Kac, V. G.: Conformal modules, Asian J.~Math., vol.~{\bf 1} no.~1 (1997), pp.~181--193.

\item{2.} D'Andrea, A.; Kac, V.~G.: Structure theory of finite conformal
algebras.

\item{3.}Fuchs, D. B.: Cohomology of infinite-dimensional Lie algebras. New York
and London: Consultants Bureau 1986.

\item{4.} Kac, V.~G.: Vertex algebras for beginners. Providence: AMS University
lecture notes vol.~{\bf 10} 1996.

\item{5.} Kac, V.~G.: The idea of locality. In: H.-D.~Doebner et al.~(eds.) Physical applications and mathematical aspects of geometry, groups and algebras, Singapore: World Scientific, 1997 pp.~16--32.

\end